\documentclass[12pt]{amsart}
\usepackage{amsfonts}
\usepackage{epsf}
\usepackage{epic}
\usepackage{eepic}
\usepackage{latexsym}
\usepackage{pstricks}

\catcode`\@=11
\long\def\@makecaption#1#2{%
    \vskip 10pt

\setbox\@tempboxa\hbox{%\ifvoid\tinybox\else\box\tinybox\fi
      \small\sf{\bfcaptionfont #1. }\ignorespaces #2}%
    \ifdim \wd\@tempboxa >\captionwidth {%
        \rightskip=\@captionmargin\leftskip=\@captionmargin
        \unhbox\@tempboxa\par}%
      \else
        \hbox to\hsize{\hfil\box\@tempboxa\hfil}%
    \fi
}
\font\bfcaptionfont=cmssbx10 scaled \magstephalf
\newdimen\@captionmargin\@captionmargin=2\parindent
\newdimen\captionwidth\captionwidth=\hsize
\catcode`\@=12

\def\lmn#1{\vadjust{\setbox1=\vtop{\hsize 12mm
\parindent=0pt\baselineskip=9pt
\rightskip=4mm plus 4mm#1}
\hbox{\kern-12mm\smash{\raise .5ex\box1}}}}

\newtheorem{theorem}{Theorem}[section]

\newtheorem{lemma}[theorem]{Lemma}
\newtheorem{proposition}[theorem]{Proposition}
\newtheorem{definition}[theorem]{Definition}
\newtheorem{corollary}[theorem]{Corollary}
\newtheorem{remark}[theorem]{Remark}

\newcommand{\eproof}{\begin{flushright} $\Box$ \end{flushright}}

\def\A{\mathcal A}

\def\Sy{\mathfrak S}
\def\P{\mathcal P}
\def\D{\mathcal D}

\def\s{\sigma}
\def\I{\mathbf I}
\def\Y{\mathrm Y}
\def\LP{{\Lambda^{(p)}}}

\def\Z{{\mathbb Z}}
\def\R{{\mathbb R}}
\def\Q{{\mathbb Q}}
\def\L{{\mathcal L}}

\def\build#1_#2^#3{\mathrel{\mathop{\kern 0pt#1}\limits_{#2}^{#3}}}
\def\mapright#1{\smash{\mathop{\longrightarrow}\limits^{#1}}}

\begin{document}

\newcommand{\MILN}{
%\begin{pspicture}[.4](-.1,-.1)(1.3,1.1)
\begin{pspicture}[0](0,.1)(.8,.6)
%\psset{unit=.5cm}
%\psgrid
\psset{unit=.5}
\psline[linewidth=.5pt,linestyle=dotted,dotsep=1pt](.1,.1)(.6,.5)
\psline[linewidth=.5pt,linestyle=dotted,dotsep=1pt](.1,.9)(.6,.5)
\psline[linewidth=.5pt,linestyle=dotted,dotsep=1pt](.6,.5)(1.1,.5)
\end{pspicture}
}

\newcommand{\MUa}{
\begin{pspicture}[.4](-.5,-.5)(1.5,1.5)
%\psgrid
%\psset{unit=1cm}
\psline[linewidth=.5pt,linestyle=dotted,dotsep=1pt](0,0)(0,1)
\psecurve[linewidth=.5pt,linestyle=dotted,dotsep=1pt](-.5,.5) (0,.5) (.3,.5)(1,0) (1.1,-.5)

\put(-.4,.9){$\scriptstyle{i}$}
\put(-.4,0){$\scriptstyle{j}$}
\put(1.2,.9){$\scriptstyle{l}$}
\put(1.2,0){$\scriptstyle{k}$}
\psdots[dotscale=.6](0,0)(0,1)(1,0)(1,1)
\end{pspicture}
}

\newcommand{\MUb}{
\begin{pspicture}[.4](-.5,-.5)(1.5,1.5)
%\psgrid
%\psset{unit=1cm}
\psline[linewidth=.5pt,linestyle=dotted,dotsep=1pt](0,0)(0,1)
\psecurve[linewidth=.5pt,linestyle=dotted,dotsep=1pt](-.5,.5) (0,.5) (.3,.5) (1,1) (1.1,1.5)

\put(-.4,.9){$\scriptstyle{i}$}
\put(-.4,0){$\scriptstyle{j}$}
\put(1.2,.9){$\scriptstyle{l}$}
\put(1.2,0){$\scriptstyle{k}$}
\psdots[dotscale=.6](0,0)(0,1)(1,0)(1,1)
\end{pspicture}
}

\newcommand{\MUc}{
\begin{pspicture}[.4](-.5,-.5)(1.5,1.5)
%\psgrid
%\psset{unit=1cm}
\psline[linewidth=.5pt,linestyle=dotted,dotsep=1pt](1,0)(1,1)
\psecurve[linewidth=.5pt,linestyle=dotted,dotsep=1pt](-.1,-.5)(0,0) (.7,.5)(1,.5)(1.5,.5)

\put(-.4,.9){$\scriptstyle{i}$}
\put(-.4,0){$\scriptstyle{j}$}
\put(1.2,.9){$\scriptstyle{l}$}
\put(1.2,0){$\scriptstyle{k}$}
\psdots[dotscale=.6](0,0)(0,1)(1,0)(1,1)
\end{pspicture}
}

\newcommand{\MUd}{
\begin{pspicture}[.4](-.5,-.5)(1.5,1.5)
%\psgrid
%\psset{unit=1cm}
\psline[linewidth=.5pt,linestyle=dotted,dotsep=1pt](1,0)(1,1)
\psecurve[linewidth=.5pt,linestyle=dotted,dotsep=1pt](-.1,1.5)(0,1) (.7,.5)(1,.5)(1.5,.5)
\put(-.4,.9){$\scriptstyle{i}$}
\put(-.4,0){$\scriptstyle{j}$}
\put(1.2,.9){$\scriptstyle{l}$}
\put(1.2,0){$\scriptstyle{k}$}
\psdots[dotscale=.6](0,0)(0,1)(1,0)(1,1)
\end{pspicture}
}

\newcommand{\FIGEX}{
\begin{pspicture}[.4](-.5,-.5)(4.5,3.5)
%\psgrid
%\psset{unit=1cm}
\pscurve[linewidth=.5pt,arrows=->,arrowscale=2](.4,1.45)(1,1.5)(1.1,1)(.7,0)(1,-1)(2,-.5)(3,-1)(3.5,.5)(2.5,1)(3.5,1)(4.5,1)(4,3.5)(2.5,3)(0,4)(-.5,2)
\psline[linewidth=.5pt,linestyle=dotted,dotsep=1pt](1,0)(2,.5)
\psline[linewidth=.5pt,linestyle=dotted,dotsep=1pt](2,.5)(3,0)
\psline[linewidth=.5pt,linestyle=dotted,dotsep=1pt](2,.5)(2,1.5)
\psline[linewidth=.5pt,linestyle=dotted,dotsep=1pt](2,1.5)(3,2)
\psline[linewidth=.5pt,linestyle=dotted,dotsep=1pt](2,1.5)(1,2)
\psline[linewidth=.5pt,linestyle=dotted,dotsep=1pt](3,2)(3.5,3)
\psline[linewidth=.5pt,linestyle=dotted,dotsep=1pt](3,2)(4,1.5)
\psline[linewidth=.5pt,linestyle=dotted,dotsep=1pt](1,2)(.5,3)
\psline[linewidth=.5pt,linestyle=dotted,dotsep=1pt](1,2)(0,1.5)
\put(-.4,1.3){$\scriptstyle{1}$}
\put(.1,2.9){$\scriptstyle{5}$}
\put(1.5,1.1){$\scriptstyle{4}$}
\put(.9,-.6){$\scriptstyle{3}$}
\put(2.9,-.6){$\scriptstyle{2}$}
\put(4.1,1.3){$\scriptstyle{7}$}
\put(3.6,2.9){$\scriptstyle{6}$}
\psdots[dotscale=.6](0,1.5)(.5,3)(2,1.5)(3.5,3)(4,1.5)(1,0)(3,0)
\end{pspicture}
}

\newcommand{\FIGEXi}{
\begin{pspicture}[.4](-.5,1)(4.5,3.5)
%\psgrid
%\psset{unit=1cm}
\psline[linewidth=.5pt,linestyle=dotted,dotsep=1pt](2,1.5)(3,2)
\psline[linewidth=.5pt,linestyle=dotted,dotsep=1pt](2,1.5)(1,2)
\psline[linewidth=.5pt,linestyle=dotted,dotsep=1pt](3,2)(3.5,3)
\psline[linewidth=.5pt,linestyle=dotted,dotsep=1pt](3,2)(4,1.5)
\psline[linewidth=.5pt,linestyle=dotted,dotsep=1pt](1,2)(.5,3)
\psline[linewidth=.5pt,linestyle=dotted,dotsep=1pt](1,2)(0,1.5)
\put(-.4,1.3){$\scriptstyle{3}$}
\put(.1,2.9){$\scriptstyle{2}$}
\put(1.8,.9){$\scriptstyle{1}$}
\put(4.1,1.3){$\scriptstyle{4}$}
\put(3.6,2.9){$\scriptstyle{5}$}
\psdots[dotscale=.6](0,1.5)(.5,3)(2,1.5)(3.5,3)(4,1.5)
\end{pspicture}
}

\newcommand{\Fi}{
\begin{pspicture}[.3](0,.2)(1,.7)
\psline[linewidth=.5pt,linestyle=dotted,dotsep=1pt](.1,.5)(.9,.5)
\end{pspicture}
}

\newcommand{\MILNijk}{
%\begin{pspicture}[.4](-.1,-.1)(1.3,1.1)
\begin{pspicture}[.4](-.3,-.3)(1.5,1.3)
%\psgrid
\psline[linewidth=.5pt,linestyle=dotted,dotsep=1pt](.1,.1)(.6,.5)
\psline[linewidth=.5pt,linestyle=dotted,dotsep=1pt](.1,.9)(.6,.5)
\psline[linewidth=.5pt,linestyle=dotted,dotsep=1pt](.6,.5)(1.1,.5)
\put(1.3,.4){$\scriptstyle{k}$}
\put(-.2,0){$\scriptstyle{j}$}
\put(-.2,.8){$\scriptstyle{i}$}
\end{pspicture}
}

\newcommand{\MILNijknew}{
%\begin{pspicture}[.4](-.1,-.1)(1.3,1.1)
\begin{pspicture}[.4](-.3,-.3)(1.5,1.3)
%\psgrid
\psline[linewidth=.5pt,linestyle=dotted,dotsep=1pt](.1,.1)(.6,.5)
\psline[linewidth=.5pt,linestyle=dotted,dotsep=1pt](.1,.9)(.6,.5)
\psline[linewidth=.5pt,linestyle=dotted,dotsep=1pt](.6,.5)(1.1,.5)
\put(1.3,.4){$\scriptstyle{k_p}$}
\put(-.2,0){$\scriptstyle{j_p}$}
\put(-.2,.9){$\scriptstyle{i_p}$}
\end{pspicture}
}

\def\KPlus{
\begin{picture}(2,2)(-1,-1)
%\put(0,0){\circle{2}}
\put(-0.707,-0.707){\vector(1,1){1.414}}
\put(0.707,-0.707){\line(-1,1){0.6}}
\put(-0.107,0.107){\vector(-1,1){0.6}}
\end{picture}}

\def\KMinus{
\begin{picture}(2,2)(-1,-1)
%\put(0,0){\circle{2}}
\put(-0.707,-0.707){\line(1,1){0.6}}
\put(0.107,0.107){\vector(1,1){0.6}}
\put(0.707,-0.707){\vector(-1,1){1.414}}
\end{picture}}

\def\KII{
\begin{picture}(2,2)(-1,-1)
%\put(0,0){\circle{2}}
\qbezier(-0.707,-0.707)(0,0)(-0.707,0.707)
\qbezier(0.707,-0.707)(0,0)(0.707,0.707)
\put(-0.607,0.607){\vector(-1,1){0.1414}}
\put(0.607,0.607){\vector(1,1){0.1414}}
\end{picture}}

\newcommand{\Wi}{
\begin{pspicture}[.4](0,0)(1,1)
\psline{->}(0,0)(0,1)
\psline{->}(1,1)(1,0)
\psline[linewidth=.5pt,linestyle=dotted,dotsep=1pt](0,.5)(1,.5)
\end{pspicture}
}

\newcommand{\Wii}{
\begin{pspicture}[.4](0,0)(1,1)
\psline{->}(0,0)(0,.4)(1,.4)(1,0)
\psline{->}(1,1)(1,.6)(0,.6)(0,1)
\end{pspicture}
}

\newcommand{\Gi}{
\begin{pspicture}[-.5](0,0)(1,.2)
\psline(0,0)(1,0)
\psdots(0,0)(1,0)
\put(1,-.3){$\scriptstyle{2}$}
\put(0,-.3){$\scriptstyle{1}$}
\end{pspicture}
}

\newcommand{\Gii}{
\begin{pspicture}[.4](0,-.2)(1,1)
\psline(0,0)(1,0)(.5,.85)(0,0)
\psdots(0,0)(1,0)(.5,.85)
\put(1.1,-.3){$\scriptstyle{2}$}
\put(-.2,-.3){$\scriptstyle{1}$}
\put(.45,1){$\scriptstyle{3}$}
\end{pspicture}
}
\newcommand{\Gia}{
\begin{pspicture}[.4](0,-.2)(1,1)
\psline(0,0)(1,0)(.5,.85)
\psdots(0,0)(1,0)(.5,.85)
\put(1.1,-.3){$\scriptstyle{2}$}
\put(-.2,-.3){$\scriptstyle{1}$}
\put(.45,1){$\scriptstyle{3}$}
\end{pspicture}
}
\newcommand{\Gib}{
\begin{pspicture}[.4](0,-.2)(1,1)
\psline(1,0)(.5,.85)(0,0)
\psdots(0,0)(1,0)(.5,.85)
\put(1.1,-.3){$\scriptstyle{2}$}
\put(-.2,-.3){$\scriptstyle{1}$}
\put(.45,1){$\scriptstyle{3}$}
\end{pspicture}
}
\newcommand{\Gic}{
\begin{pspicture}[.4](0,-.2)(1,1)
\psline(.5,.85)(0,0)(1,0)
\psdots(0,0)(1,0)(.5,.85)
\put(1.1,-.3){$\scriptstyle{2}$}
\put(-.2,-.3){$\scriptstyle{1}$}
\put(.45,1){$\scriptstyle{3}$}
\end{pspicture}
}

\newcommand{\ACE}{
\begin{pspicture}[.5](0,-.2)(3,1)
\pscircle(.5,.5){.35}
\pscircle(1.5,.5){.35}
\pscircle(2.5,.5){.35}
\psarc[arrowlength=.7,arrowsize=2pt 4]{->}(.5,.5){.35}{20}{95}
\psarc[arrowlength=.7,arrowsize=2pt 4]{->}(1.5,.5){.35}{20}{95}
\psarc[arrowlength=.7,arrowsize=2pt 4]{->}(2.5,.5){.35}{20}{95}
\psline[linewidth=.5pt,linestyle=dotted,dotsep=1pt](.85,.5)(1.15,.5)
\psline[linewidth=.5pt,linestyle=dotted,dotsep=1pt](1.85,.5)(2.15,.5)
\put(.85,.2){$\scriptstyle{1}$}
\put(1.85,.2){$\scriptstyle{2}$}
\put(2.85,.2){$\scriptstyle{3}$}
\end{pspicture}
}

\newcommand{\ACEi}{
\begin{pspicture}[.5](0,-.2)(3,1)
%\psgrid
%\psset{unit=4cm}
\psarc(.5,.5){.35}{20}{340}
\psarc(1.5,.5){.35}{20}{160}
\psarc(1.5,.5){.35}{200}{340}
\psarc(2.5,.5){.35}{200}{160}
\psarc[arrowlength=.7,arrowsize=2pt 4]{->}(1.5,.5){.35}{20}{95}
\psline(.825,.385)(1.175,.385)
\psline(1.825,.385)(2.175,.385)
\psline(.825,.615)(1.175,.615)
\psline(1.825,.615)(2.175,.615)
\end{pspicture}
}
\newcommand{\Tree}{
\begin{pspicture}[-.6](0,0)(1,.2)
\psline(0,0)(1,0)
\psdots(0,0)(.5,0)(1,0)
\put(0,-.3){$\scriptstyle{1}$}
\put(.5,-.3){$\scriptstyle{2}$}
\put(1,-.3){$\scriptstyle{3}$}
\end{pspicture}
}

\newcommand{\Rela}{
\begin{pspicture}[.4](0,0)(1,1)
\psline{->}(0,0)(0,1)
\psline{<-}(1,0)(1,1)
\psline[linewidth=.5pt,linestyle=dotted,dotsep=1pt](0,.5)(.3,.5)
\psline[linewidth=.5pt,linestyle=dotted,dotsep=1pt](.7,.5)(1,.5)
\pscircle[linewidth=.5pt,linestyle=dotted,dotsep=1pt](.5,.5){.2}
\end{pspicture}
}
\newcommand{\Relai}{
\begin{pspicture}[.4](0,0)(1,1)
\psline{->}(0,0)(0,1)
\psline{<-}(1,0)(1,1)
\end{pspicture}
}

\newcommand{\Relb}{
\begin{pspicture}[.4](0,0)(1.4,1)
\psline[linewidth=.5pt,linestyle=dotted,dotsep=1pt](0,.5)(.35,.5)
\psline[linewidth=.5pt,linestyle=dotted,dotsep=1pt](1.05,.5)(1.4,.5)
\pscircle(.7,.5){.35}
\psarc[arrowlength=.7,arrowsize=2pt 4]{->}(.7,.5){.35}{20}{95}
\end{pspicture}
}

\newcommand{\Relbi}{
\begin{pspicture}[.4](0,0)(1.4,1)
\psline[linewidth=.5pt,linestyle=dotted,dotsep=1pt](0,.5)(1.4,.5)
\end{pspicture}
}

\newcommand{\Relc}{
\begin{pspicture}[.4](0,0)(1.4,1)
\psline[linewidth=.5pt,linestyle=dotted,dotsep=1pt](0,.5)(.3,.5)
\psline[linewidth=.5pt,linestyle=dotted,dotsep=1pt](.7,.5)(1,.5)
\pscircle[linewidth=.5pt,linestyle=dotted,dotsep=1pt](.5,.5){.2}
\psline[linewidth=.5pt,linestyle=dotted,dotsep=1pt](1,.5)(1.4,.9)
\psline[linewidth=.5pt,linestyle=dotted,dotsep=1pt](1,.5)(1.4,.1)
\end{pspicture}
}

\newcommand{\Reld}{
\begin{pspicture}[.4](0,0)(1.5,1.3)
\psline[linewidth=.5pt,linestyle=dotted,dotsep=1pt](.1,.5)(.5,.5)(.5,.1)
\psline[linewidth=.5pt,linestyle=dotted,dotsep=1pt](.5,.5)(.75,.75)(1,.5)
\psline[linewidth=.5pt,linestyle=dotted,dotsep=1pt](1,.1)(1,.5)(1.4,.5)
\psline[linewidth=.5pt,linestyle=dotted,dotsep=1pt](.75,.75)(.75, 1)(1,1.3)
\psline[linewidth=.5pt,linestyle=dotted,dotsep=1pt](.75, 1)(.5,1.3)
\end{pspicture}
}

\newcommand{\Rele}{
\begin{pspicture}[.4](0,0)(1,1)
\psline[linewidth=.5pt,linestyle=dotted,dotsep=1pt](0,.1)(.3,.5)(.7,.5)(1,.1)
\psline[linewidth=.5pt,linestyle=dotted,dotsep=1pt](0,.9)(.3,.5)
\psline[linewidth=.5pt,linestyle=dotted,dotsep=1pt](.7,.5)(1,.9)
\end{pspicture}
}

\newcommand{\Relei}{
\begin{pspicture}[.4](0,.1)(1,.9)
\pscurve[linewidth=.5pt,linestyle=dotted,dotsep=1pt](0,.9)(.3,.7)(.7,.7)(1,.9)
\pscurve[linewidth=.5pt,linestyle=dotted,dotsep=1pt](0,.1)(.1,.25)(.35,.3)
\pscurve[linewidth=.5pt,linestyle=dotted,dotsep=1pt](.65,.3)(.9,.25)(1,.1)
\pscircle[linewidth=.5pt,linestyle=dotted,dotsep=1pt](.5,.3){.15}
\end{pspicture}
}

\newcommand{\Releii}{
\begin{pspicture}[.4](0,.1)(1,.9)
%\psset{unit=4cm}
%\psgrid
\pscurve[linewidth=.5pt,linestyle=dotted,dotsep=1pt](0,.9)(.4,.4)(1,.1)
\psline[linewidth=.5pt,linestyle=dotted,dotsep=1pt](0,.1)(.58,.56)
\psline[linewidth=.5pt,linestyle=dotted,dotsep=1pt](.82,.76)(1,.9)
\pscircle[linewidth=.5pt,linestyle=dotted,dotsep=1pt](.7,.66){.15}
\end{pspicture}
}

\newcommand{\Releiii}{
\begin{pspicture}[.4](0,.1)(1,.9)
%\psset{unit=4cm}
%\psgrid
\pscurve[linewidth=.5pt,linestyle=dotted,dotsep=1pt](1,.9)(.6,.4)(0,.1)
\psline[linewidth=.5pt,linestyle=dotted,dotsep=1pt](1,.1)(.42,.56)
\psline[linewidth=.5pt,linestyle=dotted,dotsep=1pt](.18,.76)(0,.9)
\pscircle[linewidth=.5pt,linestyle=dotted,dotsep=1pt](.3,.66){.15}
\end{pspicture}
}

\newcommand{\Stepa}{
\begin{pspicture}[.4](0,0)(1,1)
\pscircle(.5,.5){.35}
\psarc[arrowlength=.7,arrowsize=2pt 4]{->}(.5,.5){.35}{20}{95}
\psarc{->}(1.5,.5){.35}{120}{240}
\psline[linewidth=.5pt,linestyle=dotted,dotsep=1pt](.85,.5)(1.15,.5)
\end{pspicture}
}
\newcommand{\Stepai}{
\begin{pspicture}[.4](0,0)(1,1)
\psarc{->}(.5,.5){.35}{120}{240}
\end{pspicture}
}

\newcommand{\Stepaa}{
\begin{pspicture}[.4](0,0)(1,1)
%\psset{unit=4cm}
%\psgrid
\pscircle(.5,.5){.35}
\psarc[arrowlength=.7,arrowsize=2pt 4]{->}(.5,.5){.35}{20}{95}
\psarc{->}(1.7,.5){.35}{120}{240}
\psline[linewidth=.5pt,linestyle=dotted,dotsep=1pt](.85,.5)(1.05,.5)(1.4,.3)
\psline[linewidth=.5pt,linestyle=dotted,dotsep=1pt](1.05,.5)(1.4,.7)
\end{pspicture}
}

\newcommand{\Stepb}{
\begin{pspicture}[.4](0,0)(1.5,1)
%\psset{unit=4cm}
%\psgrid
\pscircle(.5,.5){.35}
\psarc[arrowlength=.7,arrowsize=2pt 4]{->}(.5,.5){.35}{20}{95}
\psarc{->}(1.5,.25){.35}{150}{210}
\psarc{->}(1.5,.75){.35}{150}{210}
\psline[linewidth=.5pt,linestyle=dotted,dotsep=1pt](.8,.68)(1.15,.74)
\psline[linewidth=.5pt,linestyle=dotted,dotsep=1pt](.8,.32)(1.15,.26)
\end{pspicture}
}

\newcommand{\Stepbi}{
\begin{pspicture}[.4](0,0)(1.5,1)
%\psset{unit=4cm}
%\psgrid
\psarc{->}(1.5,.25){.35}{150}{210}
\psarc{->}(1.5,.75){.35}{150}{210}
\pscurve[linewidth=.5pt,linestyle=dotted,dotsep=1pt](1.15,.74)(.9,.6)(.9,.4)(1.15,.26)
\end{pspicture}
}

\newcommand{\Stepc}{
\begin{pspicture}[.4](0,0)(1.5,1)
%\psset{unit=4cm}
%\psgrid
\pscircle(.5,.5){.35}
\psarc[arrowlength=.7,arrowsize=2pt 4]{->}(.5,.5){.35}{20}{95}
\psarc{->}(1.8,.25){.35}{150}{210}
\psarc{->}(1.8,.75){.35}{150}{210}
\psline[linewidth=.5pt,linestyle=dotted,dotsep=1pt](.85,.5)(1.15,.5)(1.45,.26)
\psline[linewidth=.5pt,linestyle=dotted,dotsep=1pt](1.15,.5)(1.45,.74)
\psarc{->}(-.7,.5){.35}{300}{60}
\psline[linewidth=.5pt,linestyle=dotted,dotsep=1pt](-.35,.5)(.15,.5)
\end{pspicture}
}
\newcommand{\Stepci}{
\begin{pspicture}[.4](0,0)(1.5,1)
\psarc{->}(1.8,.25){.35}{150}{210}
\psarc{->}(1.8,.75){.35}{150}{210}
\psline[linewidth=.5pt,linestyle=dotted,dotsep=1pt](.85,.5)(1.15,.5)(1.45,.26)
\psline[linewidth=.5pt,linestyle=dotted,dotsep=1pt](1.15,.5)(1.45,.74)
\psarc{->}(.3,.5){.35}{300}{60}
\psline[linewidth=.5pt,linestyle=dotted,dotsep=1pt](.65,.5)(1.15,.5)
\end{pspicture}
}

\newcommand{\Stepd}{
\begin{pspicture}[.4](0,0)(1.5,1)
%\psset{unit=4cm}
%\psgrid
\pscircle(.5,.5){.35}
\psarc[arrowlength=.7,arrowsize=2pt 4]{->}(.5,.5){.35}{20}{95}
\psarc{->}(1.8,.25){.35}{150}{210}
\psarc{->}(1.8,.75){.35}{150}{210}
\psline[linewidth=.5pt,linestyle=dotted,dotsep=1pt](.85,.5)(1.15,.5)(1.45,.26)
\psline[linewidth=.5pt,linestyle=dotted,dotsep=1pt](1.15,.5)(1.45,.74)
\psarc{->}(-.8,.25){.35}{330}{30}
\psarc{->}(-.8,.75){.35}{330}{30}
\psline[linewidth=.5pt,linestyle=dotted,dotsep=1pt](.15,.5)(-.15,.5)(-.45,.26)
\psline[linewidth=.5pt,linestyle=dotted,dotsep=1pt](-.15,.5)(-.45,.74)
\end{pspicture}
}

\newcommand{\Stepdi}{
\begin{pspicture}[.4](0,0)(1.5,1)
\psarc{->}(1.8,.25){.35}{150}{210}
\psarc{->}(1.8,.75){.35}{150}{210}
\psline[linewidth=.5pt,linestyle=dotted,dotsep=1pt](.85,.5)(1.15,.5)(1.45,.26)
\psline[linewidth=.5pt,linestyle=dotted,dotsep=1pt](1.15,.5)(1.45,.74)
\psarc{->}(.2,.25){.35}{330}{30}
\psarc{->}(.2,.75){.35}{330}{30}
\psline[linewidth=.5pt,linestyle=dotted,dotsep=1pt](.85,.5)(.55,.26)
\psline[linewidth=.5pt,linestyle=dotted,dotsep=1pt](.85,.5)(.55,.74)
\end{pspicture}
}

\newcommand{\Stepdii}{
\begin{pspicture}[.4](0,0)(1.5,1)
\psarc{->}(1.8,.25){.35}{150}{210}
\psarc{->}(1.8,.75){.35}{150}{210}
\psline[linewidth=.5pt,linestyle=dotted,dotsep=1pt](.55,.74)(1.45,.74)
\psarc{->}(.2,.25){.35}{330}{30}
\psarc{->}(.2,.75){.35}{330}{30}
\end{pspicture}
}

\newcommand{\Stepdiia}{
\begin{pspicture}[.4](0,0)(1.5,1)
\psarc{->}(1.8,.25){.35}{150}{210}
\psarc{->}(1.8,.75){.35}{150}{210}
\psline[linewidth=.5pt,linestyle=dotted,dotsep=1pt](.55,.26)(1.45,.26)
\psarc{->}(.2,.25){.35}{330}{30}
\psarc{->}(.2,.75){.35}{330}{30}
\end{pspicture}
}

\newcommand{\Stepdiib}{
\begin{pspicture}[.4](0,0)(1.5,1)
\psarc{->}(1.8,.25){.35}{150}{210}
\psarc{->}(1.8,.75){.35}{150}{210}
\psline[linewidth=.5pt,linestyle=dotted,dotsep=1pt](.55,.74)(1.45,.26)
\psarc{->}(.2,.25){.35}{330}{30}
\psarc{->}(.2,.75){.35}{330}{30}
\end{pspicture}
}

\newcommand{\Stepdiic}{
\begin{pspicture}[.4](0,0)(1.5,1)
\psarc{->}(1.8,.25){.35}{150}{210}
\psarc{->}(1.8,.75){.35}{150}{210}
\psline[linewidth=.5pt,linestyle=dotted,dotsep=1pt](.55,.26)(1.45,.74)
\psarc{->}(.2,.25){.35}{330}{30}
\psarc{->}(.2,.75){.35}{330}{30}
\end{pspicture}
}

\newcommand{\Treeqa}{
\begin{pspicture}[.4](-.2,-.2)(2.2,1)
%\psset{unit=4cm}
\psarc{->}(2.1,.25){.35}{150}{210}
\psarc{->}(2.1,.75){.35}{150}{210}
\psline[linewidth=.5pt,linestyle=dotted,dotsep=1pt](1.45,.5)(1.75,.26)
\psline[linewidth=.5pt,linestyle=dotted,dotsep=1pt](1.45,.5)(1.75,.74)
\psline[linewidth=.5pt,linestyle=dotted,dotsep=1pt](.55,.5)(1.45,.5)
\psarc{->}(-.1,.25){.35}{330}{30}
\psarc{->}(-.1,.75){.35}{330}{30}
\psline[linewidth=.5pt,linestyle=dotted,dotsep=1pt](.55,.5)(.25,.26)
\psline[linewidth=.5pt,linestyle=dotted,dotsep=1pt](.55,.5)(.25,.74)
\psline[linewidth=.5pt,linestyle=dotted,dotsep=1pt](1,.5)(1,.1)
\psarc{->}(1,-.25){.35}{60}{120}
\put(.65,.6){$\scriptstyle{\times}$}
\end{pspicture}
}
\newcommand{\Treeqb}{
\begin{pspicture}[.4](-.2,-.2)(2.2,1)
%\psset{unit=4cm}
\psarc{->}(2.1,.25){.35}{150}{210}
\psarc{->}(2.1,.75){.35}{150}{210}
\psline[linewidth=.5pt,linestyle=dotted,dotsep=1pt](1.45,.5)(1.75,.26)
\psline[linewidth=.5pt,linestyle=dotted,dotsep=1pt](1.45,.5)(1.75,.74)
\psline[linewidth=.5pt,linestyle=dotted,dotsep=1pt](.55,.5)(1.45,.5)
\psarc{->}(-.1,.25){.35}{330}{30}
\psarc{->}(-.1,.75){.35}{330}{30}
\psline[linewidth=.5pt,linestyle=dotted,dotsep=1pt](.55,.5)(.25,.26)
\psline[linewidth=.5pt,linestyle=dotted,dotsep=1pt](.55,.5)(.25,.74)
\psline[linewidth=.5pt,linestyle=dotted,dotsep=1pt](1,.5)(1,.1)
\psarc{->}(1,-.25){.35}{60}{120}
\put(1.15,.6){$\scriptstyle{\times}$}
\end{pspicture}
}

\newcommand{\Treeqi}{
\begin{pspicture}[.4](-.2,-.2)(2.2,1)
%\psset{unit=4cm}
\psarc{->}(2.1,.25){.35}{150}{210}
\psarc{->}(2.1,.75){.35}{150}{210}
\psline[linewidth=.5pt,linestyle=dotted,dotsep=1pt](1.45,.5)(1.75,.26)
\psline[linewidth=.5pt,linestyle=dotted,dotsep=1pt](1.45,.5)(1.75,.74)
\psarc{->}(-.1,.25){.35}{330}{30}
\psarc{->}(-.1,.75){.35}{330}{30}
\pscurve[linewidth=.5pt,linestyle=dotted,dotsep=1pt](.25,.74)(1,.5)(1.45,.5)
\psarc{->}(1,-.25){.35}{60}{120}
\end{pspicture}
}
\newcommand{\Treeqii}{
\begin{pspicture}[.4](-.2,-.2)(2.2,1)
%\psset{unit=4cm}
\psarc{->}(2.1,.25){.35}{150}{210}
\psarc{->}(2.1,.75){.35}{150}{210}
\psline[linewidth=.5pt,linestyle=dotted,dotsep=1pt](1.45,.5)(1.75,.26)
\psline[linewidth=.5pt,linestyle=dotted,dotsep=1pt](1.45,.5)(1.75,.74)
\psarc{->}(-.1,.25){.35}{330}{30}
\psarc{->}(-.1,.75){.35}{330}{30}
\pscurve[linewidth=.5pt,linestyle=dotted,dotsep=1pt](.25,.26)(1,.5)(1.45,.5)
\psarc{->}(1,-.25){.35}{60}{120}
\end{pspicture}
}

\newcommand{\Treeqiii}{
\begin{pspicture}[.4](-.2,0)(2.2,1)
%\psset{unit=4cm}
\psarc{->}(2.1,.25){.35}{150}{210}
\psarc{->}(2.1,.75){.35}{150}{210}
\psarc{->}(-.1,.25){.35}{330}{30}
\psarc{->}(-.1,.75){.35}{330}{30}
\pscurve[linewidth=.5pt,linestyle=dotted,dotsep=1pt](1.75,.74)(1,.5)(.55,.5)
\psline[linewidth=.5pt,linestyle=dotted,dotsep=1pt](.55,.5)(.25,.26)
\psline[linewidth=.5pt,linestyle=dotted,dotsep=1pt](.55,.5)(.25,.74)
\psarc{->}(1,-.25){.35}{60}{120}
\end{pspicture}
}

\newcommand{\Treeqiv}{
\begin{pspicture}[.4](-.2,0)(2.2,1)
%\psset{unit=4cm}
\psarc{->}(2.1,.25){.35}{150}{210}
\psarc{->}(2.1,.75){.35}{150}{210}
\psarc{->}(-.1,.25){.35}{330}{30}
\psarc{->}(-.1,.75){.35}{330}{30}
\pscurve[linewidth=.5pt,linestyle=dotted,dotsep=1pt](1.75,.26)(1,.5)(.55,.5)
\psline[linewidth=.5pt,linestyle=dotted,dotsep=1pt](.55,.5)(.25,.26)
\psline[linewidth=.5pt,linestyle=dotted,dotsep=1pt](.55,.5)(.25,.74)
\psarc{->}(1,-.25){.35}{60}{120}
\end{pspicture}
}

\newcommand{\treea}{
\begin{pspicture}[.5](-.2,-.2)(3.6,1)
%\psset{unit=2cm}
%\psgrid
\psarc{->}(-.1,.25){.35}{330}{30}
\psarc{->}(-.1,.75){.35}{330}{30}
\psline[linewidth=.5pt,linestyle=dotted,dotsep=1pt](.55,.5)(.25,.26)
\psline[linewidth=.5pt,linestyle=dotted,dotsep=1pt](.55,.5)(.25,.74)
\psline[linewidth=.5pt,linestyle=dotted,dotsep=1pt](.55,.5)(3.25,.5)
\psline[linewidth=.5pt,linestyle=dotted,dotsep=1pt](1,.5)(1,.1)
\psarc{->}(1,-.25){.35}{60}{120}
\put(1.15,.6){$\scriptstyle{\times}$}
\psline[linewidth=.5pt,linestyle=dotted,dotsep=1pt](1.5,.5)(1.5,.1)
\psarc{->}(1.5,-.25){.35}{60}{120}
\psline[linewidth=.5pt,linestyle=dotted,dotsep=1pt](2,.5)(2,.1)
\psarc{->}(2,-.25){.35}{60}{120}
\psline[linewidth=.5pt,linestyle=dotted,dotsep=1pt](2.8,.5)(2.8,.1)
\psarc{->}(2.8,-.25){.35}{60}{120}
\psarc{->}(3.9,.25){.35}{150}{210}
\psarc{->}(3.9,.75){.35}{150}{210}
\psline[linewidth=.5pt,linestyle=dotted,dotsep=1pt](3.25,.5)(3.55,.26)
\psline[linewidth=.5pt,linestyle=dotted,dotsep=1pt](3.25,.5)(3.55,.74)
\put(2.15,.2){$\cdots$}
\end{pspicture}
}

\newcommand{\treeb}{
\begin{pspicture}[.5](-.2,-.2)(3.6,1)
%\psset{unit=2cm}
%\psgrid
\psarc{->}(-.1,.25){.35}{330}{30}
\psarc{->}(-.1,.75){.35}{330}{30}
\psline[linewidth=.5pt,linestyle=dotted,dotsep=1pt](.55,.5)(.25,.26)
\psline[linewidth=.5pt,linestyle=dotted,dotsep=1pt](.55,.5)(.25,.74)
\psline[linewidth=.5pt,linestyle=dotted,dotsep=1pt](.55,.5)(3.25,.5)
%\psline[linewidth=.5pt,linestyle=dotted,dotsep=1pt](1,.5)(1,.1)
\psarc{->}(1,-.25){.35}{60}{120}
%\put(1.15,.6){$\scriptstyle{\times}$}
%\psline[linewidth=.5pt,linestyle=dotted,dotsep=1pt](1.5,.5)(1.5,.1)
\psarc{->}(1.5,-.25){.35}{60}{120}
\psline[linewidth=.5pt,linestyle=dotted,dotsep=1pt](2,.5)(2,.1)
\psarc{->}(2,-.25){.35}{60}{120}
\psline[linewidth=.5pt,linestyle=dotted,dotsep=1pt](2.8,.5)(2.8,.1)
\psarc{->}(2.8,-.25){.35}{60}{120}
\psarc{->}(3.9,.25){.35}{150}{210}
\psarc{->}(3.9,.75){.35}{150}{210}
\psline[linewidth=.5pt,linestyle=dotted,dotsep=1pt](3.25,.5)(3.55,.26)
\psline[linewidth=.5pt,linestyle=dotted,dotsep=1pt](3.25,.5)(3.55,.74)
\put(2.15,.2){$\cdots$}
\end{pspicture}
}

\newcommand{\Wheel}{
\begin{pspicture}[.5](-1,0)(1,1)
%\psset{unit=4cm}
%\psgrid
\pscircle[linewidth=.5pt,linestyle=dotted,dotsep=1pt](0,.5){.35}
\psline[linewidth=.5pt,linestyle=dotted,dotsep=1pt](.35,.5)(.65,.5)
\psline[linewidth=.5pt,linestyle=dotted,dotsep=1pt](-.35,.5)(-.65,.5)
\psline[linewidth=.5pt,linestyle=dotted,dotsep=1pt](.22,.77)(.5,1)
\psline[linewidth=.5pt,linestyle=dotted,dotsep=1pt](-.22,.23)(-.5,0)
\psline[linewidth=.5pt,linestyle=dotted,dotsep=1pt](.22,.23)(.5,0)
\psline[linewidth=.5pt,linestyle=dotted,dotsep=1pt](-.22,.77)(-.5,1)
\end{pspicture}
}

\newcommand{\Wheelb}{
\begin{pspicture}[.5](-1,0)(1,1)
%\psset{unit=4cm}
%\psgrid
\pscircle[linewidth=.5pt,linestyle=dotted,dotsep=1pt](0,.5){.35}
\psline[linewidth=.5pt,linestyle=dotted,dotsep=1pt](.35,.5)(.65,.5)
\psline[linewidth=.5pt,linestyle=dotted,dotsep=1pt](-.35,.5)(-.65,.5)
\psline[linewidth=.5pt,linestyle=dotted,dotsep=1pt](.22,.77)(.5,1)
\psline[linewidth=.5pt,linestyle=dotted,dotsep=1pt](-.22,.23)(-.5,0)
\psline[linewidth=.5pt,linestyle=dotted,dotsep=1pt](.22,.23)(.5,0)
\psline[linewidth=.5pt,linestyle=dotted,dotsep=1pt](-.22,.77)(-.5,1)
\psline[linewidth=.5pt,linestyle=dotted,dotsep=1pt](.65,.5)(.9,.8)
\psline[linewidth=.5pt,linestyle=dotted,dotsep=1pt](.65,.5)(.9,.2)
\end{pspicture}
}
\newcommand{\STU}{
\begin{pspicture}[.4](-.1,-.1)(3.1,1.1)
%\psgrid
\psline{->}(0,0)(0,1)
\psline[linewidth=.5pt,linestyle=dotted,dotsep=1pt](0,.5)(0.25,0.5)
\psline[linewidth=.5pt,linestyle=dotted,dotsep=1pt](0.25,0.5)(0.5,0.3)
\psline[linewidth=.5pt,linestyle=dotted,dotsep=1pt](0.25,0.5)(0.5,0.7)
\put(0.7,0.4){$=$}
\psline{->}(1.25,0)(1.25,1)
\psline[linewidth=.5pt,linestyle=dotted,dotsep=1pt](1.25,0.3)(1.75,0.3)
\psline[linewidth=.5pt,linestyle=dotted,dotsep=1pt](1.25,0.7)(1.75,0.7)
\put(2,0.4){$-$}
\psline{->}(2.5,0)(2.5,1)
\psline[linewidth=.5pt,linestyle=dotted,dotsep=1pt](2.5,0.3)(3,0.7)
\psline[linewidth=.5pt,linestyle=dotted,dotsep=1pt](2.5,0.7)(3,0.3)
\end{pspicture}
}

\newcommand{\IHX}{
\begin{pspicture}[.4](-.1,-.1)(3.1,1.1)
%\psgrid
\psline[linewidth=.5pt,linestyle=dotted,dotsep=1pt](0.25,0)(0.25,1)
\psline[linewidth=.5pt,linestyle=dotted,dotsep=1pt](0,0)(0.5,0)
\psline[linewidth=.5pt,linestyle=dotted,dotsep=1pt](0,1)(0.5,1)
\put(0.7,0.4){$=$}
\psline[linewidth=.5pt,linestyle=dotted,dotsep=1pt](1.25,0.5)(1.75,.5)
\psline[linewidth=.5pt,linestyle=dotted,dotsep=1pt](1.25,0)(1.25,1)
\psline[linewidth=.5pt,linestyle=dotted,dotsep=1pt](1.75,0)(1.75,1)
\put(2,0.4){$-$}
\psline[linewidth=.5pt,linestyle=dotted,dotsep=1pt](2.5,0)(3,1)
\psline[linewidth=.5pt,linestyle=dotted,dotsep=1pt](2.5,1)(3,0)
\psline[linewidth=.5pt,linestyle=dotted,dotsep=1pt](2.625,0.25)(2.875,0.25)
\end{pspicture}
}

\newcommand{\ASi}{
\begin{pspicture}[.4](-.1,-.1)(3.6,1.1)
%\psgrid
\psecurve[linewidth=.5pt,linestyle=dotted,dotsep=1pt](-0.1,-0.5)(0,0)(0.25,0.5)
(0.375,0.65)(0.5,0.5)(0.375,0.35)(0.25,0.5)(0,1)(-0.1,1.5)
\psline[linewidth=.5pt,linestyle=dotted,dotsep=1pt](0.5,0.5)(1,0.5)
\put(1.3,0.4){$=\ -$} 
\psline[linewidth=.5pt,linestyle=dotted,dotsep=1pt](2.5,0)(3,.5)
\psline[linewidth=.5pt,linestyle=dotted,dotsep=1pt](2.5,1)(3,.5)
\psline[linewidth=.5pt,linestyle=dotted,dotsep=1pt](3,.5)(3.5,.5)
\end{pspicture}
}

\newcommand{\ExThree}{
\begin{pspicture}[.3](-.1,-.6)(2,1.1)
%\psset{unit=4cm}
%\psgrid
\pscircle(0,.5){.35}
\psarc[arrowlength=.7,arrowsize=2pt 4]{->}(0,.5){.35}{20}{95}
\put(-.1,.56){$\scriptstyle 1$}
\pscircle(1.5,.5){.35}
\put(1.4,.56){$\scriptstyle 3$}
\psarc[arrowlength=.7,arrowsize=2pt 4]{->}(1.5,.5){.35}{20}{95}
\pscircle(.75,-.4){.35}
\put(.65,-.34){$\scriptstyle 2$}
\psarc[arrowlength=.7,arrowsize=2pt 4]{->}(.75,-.4){.35}{200}{275}
\psline[linewidth=.5pt,linestyle=dotted,dotsep=1pt](.3,.7)(1.2,.7)
\psline[linewidth=.5pt,linestyle=dotted,dotsep=1pt](.3,.3)(1.2,.3)
\psline[linewidth=.5pt,linestyle=dotted,dotsep=1pt](.6,.7)(.6,-.1)
\psline[linewidth=.5pt,linestyle=dotted,dotsep=1pt](.9,.3)(.9,-.1)
\end{pspicture}
}
\newcommand{\ExThreei}{
\begin{pspicture}[.6](-.5,-.6)(2,1.1)
%%  \psset{unit=4cm}
%% \psgrid
\pscircle(0,.5){.35}
\put(-.1,.56){$\scriptstyle 1$}
\psarc[arrowlength=.7,arrowsize=2pt 4]{->}(0,.5){.35}{20}{95}
\pscircle(1.5,.5){.35}
\put(1.4,.56){$\scriptstyle 3$}
\psarc[arrowlength=.7,arrowsize=2pt 4]{->}(1.5,.5){.35}{20}{95}
\psline[linewidth=.5pt,linestyle=dotted,dotsep=1pt](.35,.5)(1.15,.5)
\end{pspicture}
}
\newcommand{\ExThreeii}{
\begin{pspicture}[.6](-.5,-.6)(.9,1.1)
%\psset{unit=4cm}
%\psgrid
\pscircle(0,.5){.35}
\put(-.1,.56){$\scriptstyle 1$}
\psarc[arrowlength=.7,arrowsize=2pt 4]{->}(0,.5){.35}{20}{95}
\end{pspicture}
}

\newcommand{\ExFive}{
\begin{pspicture}[.6](-.5,-1.5)(2,1)
%\psset{unit=4cm}
%\psgrid
\pscircle(0,.5){.35}
\put(-.1,.56){$\scriptstyle 1$}
\psarc[arrowlength=.7,arrowsize=2pt 4]{->}(0,.5){.35}{20}{95}
\pscircle(1.5,.5){.35}
\put(1.4,.56){$\scriptstyle 5$}
\psarc[arrowlength=.7,arrowsize=2pt 4]{->}(1.5,.5){.35}{20}{95}
\pscircle(-.5,-.5){.35}
\put(-.6,-.44){$\scriptstyle 2$}
\psarc[arrowlength=.7,arrowsize=2pt 4]{->}(-.5,-.5){.35}{20}{95}
\pscircle(2,-.5){.35}
\put(1.9,-.44){$\scriptstyle 4$}
\psarc[arrowlength=.7,arrowsize=2pt 4]{->}(2,-.5){.35}{20}{95}
\pscircle(.75,-1.2){.35}
\put(.65,-1.14){$\scriptstyle 3$}
\psarc[arrowlength=.7,arrowsize=2pt 4]{->}(.75,-1.2){.35}{20}{95}
\pscurve[linewidth=.5pt,linestyle=dotted,dotsep=1pt](.35,.5)(.75,.4)(1.15,.5)
\pscurve[linewidth=.5pt,linestyle=dotted,dotsep=1pt](1.7,-.3)(1,0)(.75,.4)
\pscurve[linewidth=.5pt,linestyle=dotted,dotsep=1pt](.2,.2)(.1,-.2)(-.2,-.35)
\psline[linewidth=.5pt,linestyle=dotted,dotsep=1pt](.46,-1)(.1,-.2)
\pscurve[linewidth=.5pt,linestyle=dotted,dotsep=1pt](1.6,0.17)(1.63,0)(1.8,-0.2)
\pscurve[linewidth=.5pt,linestyle=dotted,dotsep=1pt](1.63,0)(1,-.5)(.9,-.9)
\pscurve[linewidth=.5pt,linestyle=dotted,dotsep=1pt](1.3,.2)(.6,-.4)(-.16,-.6)
\psline[linewidth=.5pt,linestyle=dotted,dotsep=1pt](.6,-.4)(.7,-.85)
\end{pspicture}
}

\newcommand{\ExFour}{
\begin{pspicture}[.6](-.5,0)(2,1)
%%\psset{unit=4cm}
%\psgrid
\pscircle(0,.5){.35}
\put(-.1,.56){$\scriptstyle 1$}
\psarc[arrowlength=.7,arrowsize=2pt 4]{->}(0,.5){.35}{20}{95}
\pscircle(1.5,.5){.35}
\put(1.4,.56){$\scriptstyle 2$}
\psarc[arrowlength=.7,arrowsize=2pt 4]{->}(1.5,.5){.35}{20}{95}
\psline[linewidth=.5pt,linestyle=dotted,dotsep=1pt](.3,.7)(1.2,.7)
\psline[linewidth=.5pt,linestyle=dotted,dotsep=1pt](.3,.3)(1.2,.3)
\psline[linewidth=.5pt,linestyle=dotted,dotsep=1pt](.75,.7)(.75,.3)
\end{pspicture}
}

\newcommand{\ExFouri}{
\begin{pspicture}[.6](-.5,-1)(2,1)
%\psset{unit=4cm}
%\psgrid
\pscircle(0,.5){.35}
\put(-.1,.56){$\scriptstyle 1$}
\psarc[arrowlength=.7,arrowsize=2pt 4]{->}(0,.5){.35}{20}{95}
\pscircle(1.5,.5){.35}
\put(1.4,.56){$\scriptstyle 4$}
\psarc[arrowlength=.7,arrowsize=2pt 4]{->}(1.5,.5){.35}{20}{95}
\psline[linewidth=.5pt,linestyle=dotted,dotsep=1pt](.3,.7)(1.2,.7)
\psline[linewidth=.5pt,linestyle=dotted,dotsep=1pt](.3,.3)(1.2,.3)
\psline[linewidth=.5pt,linestyle=dotted,dotsep=1pt](.75,.7)(.75,.3)
\pscircle(0,-.8){.35}
\put(-.1,-.74){$\scriptstyle 2$}
\psarc[arrowlength=.7,arrowsize=2pt 4]{->}(0,-.8){.35}{200}{275}
\pscircle(1.5,-.8){.35}
\put(1.4,-.74){$\scriptstyle 3$}
\psarc[arrowlength=.7,arrowsize=2pt 4]{->}(1.5,-.8){.35}{200}{275}
\psline[linewidth=.5pt,linestyle=dotted,dotsep=1pt](1.5,.15)(1.5,-.45)
\psline[linewidth=.5pt,linestyle=dotted,dotsep=1pt](0,.15)(0,-.45)
\pscurve[linewidth=.5pt,linestyle=dotted,dotsep=1pt](0,-.1)(1,-.3)(1.27,-.53)
\pscurve[linewidth=.5pt,linestyle=dotted,dotsep=1pt](1.5,-.1)(.5,-.3)(.23,-.53)
\end{pspicture}
}

\title[Milnor numbers, Trees, and Conway polynomial]
{Milnor numbers, Spanning Trees, and  the Alexander-Conway Polynomial}
\thanks{2000 \emph{Mathematics Subject Classification:} 57M27.}

\author{Gregor Masbaum}
\address{Institut de Math{\'e}matiques de Jussieu, Equipe `Topologie
et g{\'e}om{\'e}trie  alg{\'e}briques', Case 7012, Universit{\'e}
Paris VII,   75251 Paris Cedex 05, France}
\email{masbaum@math.jussieu.fr}

\author{Arkady Vaintrob}
\address{Department of Mathematics, University of Oregon, 
Eugene, OR 97405, USA} 
\email{vaintrob@math.uoregon.edu}

\date{October 2001}

\begin{abstract}  
We study relations between the Alexander-Conway 
polynomial $\nabla_L$ and Milnor higher linking numbers of 
links from the point of view of finite-type (Vassiliev) invariants. 
We give a formula for the first non-vanishing coefficient   
of $\nabla_L$ of an $m$-component link $L$ all of whose 
Milnor numbers $\mu_{i_1\ldots i_p}$ vanish for $p\le n$.  
We express this coefficient 
        as a polynomial in Milnor numbers of $L$.
Depending on whether the parity of $n$ is odd or even,
the terms in this polynomial correspond 
either to spanning trees 
in certain graphs or to decompositions of certain $3$-graphs 
into pairs of spanning trees.
  Our results complement  determinantal formulas
of Traldi and Levine obtained by geometric methods.
\end{abstract}

\maketitle

\tableofcontents

\section{Introduction}

The Alexander-Conway polynomial 
$$\nabla_L(z)= \sum_{i\geq 0} c_i(L) z^i \in \Z[z]$$ 
of a link $L$ in $\R^3$
is one of the most thoroughly studied classical isotopy invariants of 
links. 
In this paper we study relations between $\nabla_L$ and 
and the Milnor higher linking numbers of $L$
from the point of view of 
the
theory of finite-type (Vassiliev) invariants.   

Hosokawa~\cite{Hw}, Hartley~\cite[(4.7)]{Ha},  
and Hoste~\cite{Ho} showed that the coefficients $c_i(L)$ of
$\nabla_L$ for an $m$-component link $L$ vanish when $i\leq m-2$ 
and that the coefficient $c_{m-1}(L)$ depends only
on the linking numbers $\ell_{ij}(L)$ between the $i$th and $j$th 
components of $L$. Namely,
\begin{equation}
  \label{eq:ho_det}
c_{m-1}(L)= \det\Lambda^{(p)},
\end{equation}
where 
$\Lambda=(\lambda_{ij})$ is the matrix formed by linking numbers
$$\lambda_{ij}=
\begin{cases}
\ \ \ \ \ -\ell_{ij}(L), & \text{\ if \ $i\ne j$\ }  \\
\sum_{k\neq i} \ell_{ik}(L), & \text{\ if \ $i=j$}
\end{cases}
$$
and $\Lambda^{(p)}$ denotes the matrix obtained by removing from
$\Lambda$ the $p$th row and column (it is easy to see that \ 
$\det\Lambda^{(p)}$ \ does not depend on $p$).

If the link $L$ is \emph{algebraically split}, \emph{i.e.}\ all linking
numbers  $\ell_{ij}$ vanish, then not only  $c_{m-1}(L)=0$, but, as
was proved by Traldi~\cite{Tr1,Tr2} and Levine~\cite{Le1}, the next $m-2$ 
coefficients of $\nabla_L$ also vanish
$$c_{m-1}(L)=c_{m}(L)=\ldots=c_{2m-3}(L)=0.$$    

For algebraically split oriented links,
there exist well-defined integer-valued isotopy invariants
$\mu_{ijk}(L)$
called the \emph{Milnor triple linking numbers}.
These invariants  generalize  ordinary linking numbers
(see Section~\ref{sec:milnor}), 
but unlike $\ell_{ij}$, the triple linking numbers  are antisymmetric
with respect to their indices,
$\mu_{ijk}(L)=-\mu_{jik}(L)=\mu_{jki}(L).$  
Thus, for an algebraically split link $L$ with $m$ components,   
we have $m \choose 3$ triple linking numbers $\mu_{ijk}(L)$
corresponding  to the different $3$-component sublinks of $L$.

Levine~\cite{Le1} (see also Traldi~\cite[Theorem~8.2]{Tr2}) found an 
expression for the 
coefficient $c_{2m-2}(L)$ of $\nabla_L$ for an algebraically split 
$m$-component  link in terms  of triple Milnor numbers
\begin{equation}
  \label{eq:lev}
c_{2m-2}(L)=\det\LP,
\end{equation}
where $\Lambda=(\lambda_{ij})$ is an $m\times m$ skew-symmetric matrix
with entries
$$\lambda_{ij}=\sum_{k} \mu_{ijk}(L),$$
and $\Lambda^{(p)}$, as before, is the result of removing the $p$th 
row and column. 
(In particular, since the determinant of an antisymmetric matrix of
odd size is always zero, this implies that $c_{2m-2}(L)=0$   
if $m$ is  even.\footnote{This follows, of course, already from the
well-known fact 
that the Alexander-Conway polynomial of an $m$-component
link has non-zero terms only in degrees congruent to $m-1$ 
{\em mod} $2$.})
\medskip

It is well known that the coefficient $c_n$ of $\nabla(z)$ is
a finite-type (Vassiliev) invariant of order $n$ whose weight system  
can be computed by a simple combinatorial rule (see
Section~\ref{sec:ac}). 
This leads to a simple proof of
vanishing of $c_k(L)$ for $k\le m-2$  and 
also to an explicit expression for the 
coefficient $c_{m-1}$ as a sum of
products of linking numbers corresponding to
maximal trees in the complete graph $K_m$ 
with vertices $\{1,\ldots,m\}$.
Thus, there is a
second formula  for $c_{m-1}$ 
(see~\cite{Ha} and~\cite{Ho})
\begin{equation}
\label{eq:ho_tree}
c_{m-1}(L)=\D_m(\ell_{ij}(L))~.
\end{equation}
Here
$\D_m(x_{ij})$
is the \emph{Kirchhoff polynomial} 
whose monomials correspond to spanning trees in
the complete graph $K_m$ (see Section~\ref{sec:statmt}). 

The equality of the expressions~(\ref{eq:ho_det})
and~(\ref{eq:ho_tree}) for  $c_{m-1}(L)$ follows from the classical
\emph{Matrix-Tree Theorem} (see \emph{e.g.}~\cite{Bol,Tutte}) applied
to the graph $K_m$.
\medskip

Formula~(\ref{eq:lev}) is 
similar to the first determinantal expression~(\ref{eq:ho_det}).  
One of our initial goals in this work
was to find an analog of the tree sum formula~(\ref{eq:ho_tree}) for
algebraically split links,  and one of our results is a formula
expressing $c_{2m-2}$  as  the square of  a sum over trees
\begin{equation} \label{eq:pmsquare}
c_{2m-2}=\left(\P_m(\mu_{ijk})\right)^2,
\end{equation}
where $\P_m$ is the Pfaffian-tree polynomial 
 introduced in~\cite{MV}.
Similarly to the Kirchhoff polynomial $\D_m$,
the polynomial $\P_m$ is the generating function of spanning trees in the
 complete $3$-graph $\Gamma_m$ with $m$ vertices (see
 Section~\ref{sec:statmt}).

Formula~(\ref{eq:pmsquare}) can be deduced from~(\ref{eq:lev}) by our
Pfaffian Matrix-Tree Theorem proved in~\cite{MV}. Here we give a 
direct proof based on the theory of finite type invariants.

For $m=3$, both formula~\eqref{eq:lev} and our theorem give the same  
known result first proved by Cochran\footnote{Before that, Kidwell and
Morton showed that the coefficient  $c_{m-2}$  for algebraically split
links is a perfect square.}  in~\cite[Theorem 5.1]{Co} 
$$
c_4(L)=(\mu_{123}(L))^2~.
$$

In the $m\ge 5$ case, our formula is new. For example, when $m=5$, 
we obtain that the first non-vanishing coefficient of $\nabla_L(z)$ 
for algebraically split links with $5$ components is equal to 
\begin{eqnarray*} 
&c_8(L)\hspace{-6pt}&= \P_5(\mu_{ijk}(L))^2\\
&&=
\bigl( \mu_{123}(L)\mu_{145}(L) - \mu_{124}(L)\mu_{135}(L)
+ \mu_{125}(L)\mu_{134}(L)  \ \pm \ \ldots\bigr)^2, 
\end{eqnarray*}
where $\P_5(\mu_{ijk}(L))$ consists of $15$ terms corresponding to
the spanning trees in the complete $3$-graph with $5$ vertices. 

It is worth pointing out that  our  formula~(\ref{eq:pmsquare})
 has some computational advantages over~\eqref{eq:lev}.
 For example, the straightforward expansion of 
$\det\Lambda^{(p)}$ in~(\ref{eq:lev}) for $m=5$ 
would consist (before cancellations) of $729$ terms,
 each of which is a product of
four $\mu_{ijk}$'s. In contrast, the computation  based 
on~(\ref{eq:pmsquare}) only requires finding the sum of $15$ products  
of two  $\mu_{ijk}$'s each,  and taking the square.    
\medskip

If all triple Milnor linking numbers of a link $L$ vanish, then  
its Milnor numbers of length four 
are well-defined, \emph{etc.}, and so it is natural
to ask whether the lowest term of the Alexander-Conway polynomial
$\nabla_L(z)$ can be expressed via the first non-vanishing Milnor
numbers of $L$. 
If $L$ is an $m$-component link 
all of whose Milnor higher linking   
numbers $\mu_{i_1\ldots i_p}$ vanish for $p\le n$,
then, as it was proved by Traldi~\cite{Tr2}
and Levine~\cite{Le2},
the polynomial $\nabla_L(z)$ begins with the
term of degree  $n(m-1)$.
The corresponding coefficient can be computed,  as before,  by the
determinantal formula  
\begin{equation} \label{eq:detgeneral}
c_{n(m-1)}(L)=\det\Lambda^{(p)},
\end{equation}
 where the $(i,j)$th entry of the matrix $\Lambda$ is expressed via
Milnor invariants of length $n+1$
$$
\lambda_{ij}=\sum_{r_1,\ldots,r_{n-1}} \mu_{r_1,\ldots,r_{n-1},j,i}(L)~. 
$$ 

Our approach based on the theory of finite-type invariants leads to
a new proof of vanishing of the coefficients
$c_i(L)$ for $i<n(m-1)$ 
of the Alexan\-der-Conway polynomial of such a link 
and to a tree-sum counterpart of formula~(\ref{eq:detgeneral}).
Namely, we show here that  $c_{n(m-1)}(L)$ 
is equal either to the Kirchhoff polynomial or to the square of the
Pfaffian-tree polynomial whose variables are linear combinations of
the first non-vanishing Milnor numbers of $L$. 
Our results also lead to a similar formula for the next coefficient 
$c_{n(m-1)+1}(L)$ in the case 
when
both $n$ and $m$ are even (and,
therefore, $c_{n(m-1)}(L)$ is always zero); this formula is
described in Section~\ref{sec:even}.
\medskip

Whereas the original proofs of
the determinantal formulas~(\ref{eq:ho_det}),~(\ref{eq:lev})
and~(\ref{eq:detgeneral}) use classical topological 
tools (Seifert surfaces, geometric interpretation of the Milnor
linking numbers, \emph{etc.}), our approach based on 
finite-type invariants is, in some sense, more of a combinatorial 
nature. 

The connection between the Alexander-Conway polynomial and the Milnor
numbers is established by studying their weight systems and then using
the Kontsevich integral. In the dual language of the space of chord
diagrams, the Milnor numbers correspond to the tree diagrams
(see~\cite{HM}) and the Alexander-Conway polynomial can be described
in terms of certain trees and wheel diagrams (see~\cite{KSA}
and~\cite{Vai}). From this point of view, 
the trees in our formulas for the first non-vanishing term of
$\nabla_L$ appear very naturally 
and the corresponding (less symmetric) determinantal expressions
follow by the Matrix-Tree Theorems. (In fact, it was an attempt to
find a diagrammatic interpretation of formula~\eqref{eq:lev}
that led us to the discovery of the Pfaffian Matrix-Tree Theorem
for $3$-graphs~\cite{MV}.)   

The work presented in this paper is a first step in the study of 
the Alexander-Conway polynomial of links by means of finite
type invariants. 
We believe that our approach will lead to a
diagrammatic  (and perhaps more canonical) 
interpretation of  a factorization of this polynomial given by 
Levine in~\cite{Le2}.  
The factors are the Alexander-Conway polynomial of a knot obtained by
banding together the components of the link, and 
a power series involving the Milnor invariants of a certain string
link representative of the link. They should correspond respectively
to the wheel and tree parts of the weight system of $\nabla_L$.  

\medskip
The paper is organized as follows.
In Section~\ref{sec:statmt} we define the spanning-tree polynomials
$\D_m$ and $\P_m$ and formulate our results.
In Section~\ref{sec:ac}  we recall
some facts about the Alexander-Conway polynomial $\nabla_L(z)$ and its
weight system and show, as a warm-up, how the philosophy of finite-type
invariants leads naturally to 
the tree-sum formula~\eqref{eq:ho_tree} for 
$c_{m-1}(L)$. 
In Section~\ref{sec.main} we study the 
weight systems corresponding to the coefficients of $\nabla_L(z)$ and   
prove that they vanish on a certain class of diagrams. 
The proof of our Vanishing Lemma~\ref{2.1} uses properties of the
Alexander-Conway weight system 
from~\cite{FKV} which are based on the connection between $\nabla$ and
the Lie superalgebra  $gl(1|1)$. 
Milnor linking numbers and their connection with finite-type
invariants and the Kontsevich integral found in~\cite{HM}
are described in Section~\ref{sec:milnor}. 
In Section~\ref{sec:coeff} using the Vanishing Lemma of
Section~\ref{sec.main} and the methods of~\cite{HM}, we show that
for a link $L$ whose Milnor invariants of degree $\le n-1$ vanish,
the first non-vanishing coefficient of $\nabla_L(z)$ 
can be expressed as a polynomial $F_m^{(n)}$ in the 
degree-$n$ Milnor numbers of $L$.
In Section~\ref{sec:fm} we  express
this polynomial in terms of the spanning tree polynomials   
$\D_m$ and $\P_m$ introduced in 
Section~\ref{sec:statmt} and complete the proofs of
Theorems~\ref{thm:alex-trees} and~\ref{thm:general}.    
This identification is based on some algebraic properties of the
Pfaffian-tree polynomial $\P_m$ which have been established in~\cite{MV}.
In particular, one has $F_m^{(2)}=\P_m^2$, and we explain how this allows one
to compute the coefficients of
              the polynomial $F_m^{(2)}$
by counting tree decompositions of certain associated $3$-graphs.
Finally in Section~\ref{sec:even} we give a formula for
the coefficient $c_{n(m-1)+1}$ in the case where  
$c_{n(m-1)}$ is identically zero for parity reasons. 
\medskip

\noindent 
\emph{Acknowledgements.}
Much of the work on this paper has been done during the visits of
G.M.\ to the University of Oregon  and of A.V.\ 
to the University of Paris VII and IHES. We are grateful to 
these institutions for hospitality and support.
Research of the second author was partially supported by NSF grant
DMS-0104397.

\section{Spanning tree polynomials and statement of results}
\label{sec:statmt}

Recall the  following result of Hartley~\cite{Ha} and Hoste~\cite{Ho}.

\begin{theorem}\label{Hoste.thm} Let $L$ be an
oriented link in  $S^3$ with $m$ (numbered) components. Let
$\ell_{ij}(L)$ be the linking number 
between the $i$th and $j$th components. Then the 
first $m-1$ coefficients of the Alexander-Conway polynomial vanish,  
$$c_i(L)=0 \text{\ for \ }  i\leq m-2,$$
and the coefficient $c_{m-1}(L)$
is equal to the 
Kirchhoff polynomial~{\em (\ref{eq:kirch_pol})}
evaluated at the linking numbers of $L$:    
\begin{equation} \label{sumtrees.Hoste} 
c_{m-1}(L)= \D_m(\ell_{ij}(L))~.
\end{equation} 
\end{theorem}

\medskip

Here the \emph{Kirchhoff polynomial} (or the \emph{spanning trees 
generating function}) of the complete graph
$K_m$ with $m$ vertices is the polynomial
\begin{equation}
   \label{eq:kirch_pol}
\mathcal{D}_m(x_{ij})= \sum_T x_T
 \end{equation}
whose $m \choose 2$ variables 
$$x_{ij}, \ 1 \le i,j \le m, \ i \ne j, \ x_{ij}=x_{ji},$$  
correspond to the edges of 
      $K_m$ and whose
monomials $$ x_T=\prod_{e\in T}x_e$$
correspond to the maximal (spanning) trees $T$ in $K_m$. 

\begin{figure}[h]
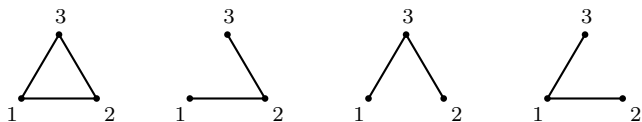

\begin{center}
\Gii \ \ \ \ \ \ \ \Gia\ \ \ \ \ \ \ \ \Gib\ \ \ \ \ \ \ \ \Gic
\caption{\label{Gamma3}The complete graph $K_3$ and its three spanning
trees.} 
\end{center}
\end{figure}

For example, if $m=2$, then $\D_2=x_{12}$,
corresponding to 
the only spanning tree in 
 $$K_2=\ \ \ \Gi$$
and so $c_1(L)=\ell_{12}(L)$. 
If $m=3$, then 
$\D_3=x_{12}x_{23} +x_{23}x_{31} +x_{31}x_{12}$ 
(see Figure \ref{Gamma3}) and  so
$$c_2(L)=\ell_{12}(L)\ell_{23}(L) +\ell_{23}(L)\ell_{31}(L)
+\ell_{31}(L)\ell_{12}(L)~.$$ 

\medskip

To state an analog of Theorem~\ref{Hoste.thm} for algebraically
split links we need to introduce another tree-generating polynomial
analogous to the Kirchhoff polynomial.

Namely, instead of usual graphs whose edges can be thought of as
 segments  joining pairs of points, we consider \emph{$3$-graphs}  
 whose edges have three (distinct) vertices and can be visualized as
 triangles or Y-shaped objects \MILN with the three vertices at their
 endpoints. 

Similarly to variables $x_{ij}$ of $\D_m$, for each triple of distinct
numbers $i,j,k \in \{1,2,\ldots,m\}$ we introduce   variables
$y_{ijk}$ antisymmetric in $i,j,k$ 
$$
  y_{ijk}=-   y_{jik}=   y_{jki}, \text{\ and\ } y_{iij}=0~.
$$
These variables correspond to edges $\{i,j,k\}$ of the \emph{complete
  $3$-graph} $\Gamma_m$ with vertices $\{1,\ldots,m\}$. 
As in the case of ordinary graphs, the correspondence
\begin{equation*}
\text{variable \ }    y_{ijk} \quad  \mapsto \quad  \text{edge\ }
\{i,j,k\} \  \text{of \ } \Gamma_m
\end{equation*}
assigns to each monomial in $  y_{ijk}$  a  sub-$3$-graph
of $\Gamma_m$. However, because of the antisymmetry, the
correspondence between monomials and sub-$3$-graphs
is not one-to-one. A sub-$3$-graph
determines a monomial only up to sign, and to define an
analog of the tree-generating function for the complete $3$-graph we
need to fix signs. 
\medskip

The generating function of spanning trees in the complete $3$-graph
with $m$  vertices, or the \emph{Pfaffian-tree polynomial} $\P_m$ is
defined as follows~\cite{MV}.  
If $m$ is even, then we set  
$$
\P_m=0
$$
(there are no spanning trees in $3$-graphs with even number of vertices). 
If $m$ is odd, then 
\begin{equation}
  \label{eq:pm}
\P_m=\sum_T \varepsilon(T)  y_T~, 
\end{equation}
where the sum is taken over all monomials
$$
 y_T=\prod_{p=1}^{d}  y_{i_p j_p k_p}
$$
of degree $d=(m-1)/2$ in $ y_{ijk}$ 
(with the convention that monomials
which differ only by changing the order of indices in some of the
variables  are taken only once) and the 
coefficient
$\varepsilon(T)\in\{0,\pm1\}$ 
for the collection of triples 
$$T=((i_1 j_1 k_1), \ldots, (i_d j_d k_d))$$
is defined by the following rule.

\begin{figure}[h]
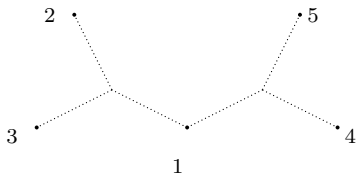

\begin{center}
\FIGEXi
\caption{\label{figexi} A spanning tree in the complete $3$-graph
$\Gamma_5$. It   has two edges, $\{1,2,3\}$ and $\{1,4,5\}$, and
contributes the term  $ y_{123}\, y_{145}$ to $\P_5$.}  
\end{center}
\end{figure}

Glue together the \ $d$ \  Y-shaped objects   
$$ \MILNijknew $$
corresponding to the variables
in the monomial \ $ y_T$ \ (see \emph{e.g.}\ Figure~\ref{figexi}).
If the resulting $1$-complex is a tree (\emph{i.e.}\ it is connected
 and  simply connected, or, equivalently, the corresponding 
sub-$3$-graph  in  $\Gamma_m$ is a tree), then 
the product
$$\sigma_T=\sigma_1\ldots \sigma_d$$   
obtained by multiplying (in an arbitrary order)
the $3$-cycles $\sigma_p=(i_p j_p k_p)$ in the
symmetric group $\Sy_m$ is an $m$-cycle 
$$\sigma=(s(1)s(2)\ldots s(m)), $$
for some $s\in \Sy_m$
and we define
\begin{equation} \label{eq:signs}
\varepsilon(T)=(-1)^s~.
\end{equation}
If the monomial $T$ does not produce a  tree,
we set $\varepsilon(T)=0$.

\begin{remark}\rm
In fact, $\sigma_T$ is an $m$-cycle 
if and only if the monomial $y_T$ corresponds to a tree
(see~\cite{MV}). 
The terms of the Kirchhoff polynomial $\D_m$ also admit a similar
description: a 
monomial $x_{i_1 j_1}\ldots x_{i_{m-1} j_{m-1}}$ enters $\D_m$ 
(always with the plus sign) if and only if the 
product of the $m-1$ transpositions $\tau_k=(i_k j_k)$ (taken in any 
order) in the symmetric group $\Sy_m$ is an $m$-cycle. 
\end{remark}

For example, if $m=5$, we have
\begin{equation} \label{m=5}
\P_5=   y_{123}\, y_{145} -  y_{124}\, y_{135} +  y_{125}\, y_{134}
\  \pm \ \ldots~,
\end{equation}
where the right-hand side is a sum of $15$ similar terms corresponding
to the $15$ spanning trees of $\Gamma_5$.

If we visualize the edges of $\Gamma_m$ as  Y-shaped
objects \MILN, then the spanning tree corresponding to the first term 
of~(\ref{m=5})  
will look like on Figure~\ref{figexi}.  

In this paper we prove the following
two theorems analogous to Theorem~\ref{Hoste.thm} for algebraically
split links. 
\begin{theorem}  \label{thm:alex-trees} 
Let $L$ be an algebraically split oriented link with
$m$ components. Then 
  \begin{equation}
  \label{eq:sqtrees}
  c_{2m-2}(L)= \bigl(\P_m(\mu_{ijk}(L))\bigr)^2~.
  \end{equation}
\end{theorem}

In fact, it follows from our \emph{Pfaffian Matrix-Tree
Theorem}~\cite{MV} (which is an analog for $3$-graphs
of the classical Matrix-Tree Theorem~\cite{Bol,Tutte}) 
that  the polynomial $\P_m(\mu_{ijk}(L))$ is equal to the Pfaffian
of the  skew-symmetric matrix $\Lambda^{(p)}$.
This gives a combinatorial explanation of the equality of the
expressions~\eqref{eq:lev} and~\eqref{eq:sqtrees}.
However, we will not use the Pfaffian Matrix-Tree Theorem
in the present paper and will give a direct proof of
Theorem~\ref{thm:alex-trees} based on the theory of finite type
invariants. 
\medskip

\begin{theorem}\label{thm:general}
Let $L$ be an oriented $m$-component link 
with vanishing Milnor numbers of length $p \le n$
and let 
$$\nabla_L(z)=\sum_{i\geq 0} c_i(L) z^i $$
be its Alexander-Conway polynomial. Then
$c_i=0$ for $i<n(m-1)$ 
and
\begin{equation}\label{eq:genform}
c_{n(m-1)}(L)= 
\begin{cases} \ \ \D_m(\ell^{(n)}_{ij})~,   &\text{if $n$ is odd}\\                       (\P_m(\mu^{(n)}_{ijk}))^2\, , &\text{if $n$ is even},                 \end{cases}
\end{equation}
where $\ell^{(n)}_{ij}$ and $\mu^{(n)}_{ijk}$ are certain linear
combinations of the Milnor numbers of $L$ of length $n+1$.
\end{theorem}

In the case when both $n$ and $m$ are even 
the coefficient $c_{n(m-1)}(L)$ is always zero
for degree reasons and so $c_{n(m-1)+1}(L)$
becomes the first non-vanishing coefficient of $\nabla_L$.
In Section~\ref{sec:even} we give a formula for
this coefficient analogous to~(\ref{eq:genform}).

\section{Alexander-Conway polynomial and spanning trees} 
\label{sec:ac}

In this section we recall the basic properties of the Alexander-Conway
polynomial $\nabla_L$ and its weight system and explain the
spanning-trees formula~(\ref{eq:ho_tree}) for the first non-vanishing
term of $\nabla_L$ from the point of view of finite-type invariants. 

\begin{figure}[h]
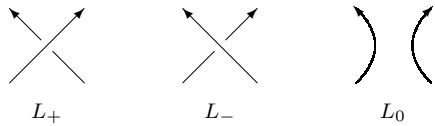

\begin{center}
\begin{displaymath}
\mathop{\KPlus}_{L_+}\qquad
\mathop{\KMinus}_{L_-}\qquad
\mathop{\KII}_{L_0}
\end{displaymath}
\caption{\label{skeintriple} A skein triple.}
\end{center}
\end{figure}

Let 
$$\nabla_L(z)=\sum_{i\geq 0} c_i(L) z^i $$
be 
the Alexander-Conway polynomial of an
oriented link $L$  in $S^3$.   It satisfies  the {\em skein
  relation} 
\begin{equation}
\nabla_{L_+}-\nabla_{L_-}=z \nabla_{L_0}~,\label{skein}
\end{equation} 
where $(L_+, L_-, L_0)$ is any skein triple (see
Figure~\ref{skeintriple}).  

The Alexander-Conway polynomial is uniquely determined by
relation~(\ref{skein}) and the initial conditions
\begin{equation}
\nabla_{U_m}=\begin{cases} 1 \ \ \ \text{if $m=1$}\\
0\ \ \ \text{if $m\geq 2$},
\end{cases}
\label{unlink}
\end{equation}
where $U_m$ is the trivial link with $m$ components.

The skein relations~(\ref{skein}) and~(\ref{unlink})
immediately give the well-known fact that
for an $m$-component link
\begin{equation} \label{eq:parvanish}
c_i(L)=0 \ \ \text{if  \ $i\equiv m$  \emph{mod}  $2$}.  
\end{equation}
They also easily imply
that for a $2$-component link $L$, one has  
$$
\nabla_L = \ell_{12}(L)z + \ldots~,
$$
where $\ell_{12}(L)$ is the linking number between the components of 
$L$. 

Theorem~\ref{Hoste.thm} generalizes this fact to links
with arbitrary number of components. One of the starting points 
for the present paper was the observation that the appearance of
spanning trees in this theorem 
is very natural from the point of view of the theory of finite type 
invariants.  As a warm-up, and also to motivate the techniques that we
use in the subsequent sections, let us briefly discuss how one can
understand  this using the properties of the Alexander-Conway weight 
system.

It is well known that the coefficient $c_n$ of the Alexander-Conway
polynomial is a finite type invariant of degree $n$. Therefore it has
a {\em   weight system\/} $W_n$ which is a linear form on the vector
space of chord diagrams of degree $n$ on $m$ circles, where $m$ is the
number of components of the link. 
Here, by a chord diagram of degree $n$ 
we understand a collection
of $n$ pairs  of points on the union of these $m$ circles. It is
represented pictorially by $n$ {\em dashed} lines, each line
connecting the two points in a pair.  

The skein relation~(\ref{skein}) implies the following 
formula 
\begin{equation}
W_{n}(\ \ \Wi\ \ ) = W_{{n-1}}(\ \ \Wii\ \ ) \label{smooth}
\end{equation} 
which allows to compute these weight systems recursively. 
For
a chord diagram $D$ on $m$ circles with $n$        
chords, let $D'$ be the result of smoothing 
of all chords by means of~(\ref{smooth}). 
Note that $D'$ is a disjoint union of several, say 
$m'$ 
circles (no chords are left). Thus,

\begin{equation}W_{n}(D)=W_{0}(D')=\begin{cases} 1 \ \ \ \text{if
$m'=1$ }\\ 
0\ \ \ \text{if 
$m'\geq 2$},
\end{cases}
\label{compW}
\end{equation}
where the second equality follows from~(\ref{unlink}).

To see how this relates to 
formula~(\ref{eq:ho_tree}),  let us compute 
$W_{n}(D)$ for a chord diagram $D$ of degree $n$ on $m$ circles. 
From~(\ref{compW}) 
we see
that $W_{n}(D)$ can only be non-zero if the diagram $D'$ obtained by
applying~(\ref{smooth}) to
all the chords of $D$ consists of just one circle.  
Since a smoothing of a chord cannot reduce the number of
circles by more than one, this means that we need at least $m-1$ chords. 
Moreover, the diagrams $D$ with exactly $m-1$ chords 
satisfying $W_{m-1}(D)\ne 0$
must have the property that if each circle
of $D$ is shrinked to a  point, the resulting graph formed by the
chords\footnote{%
In general, 
this  will actually be a multi-graph, that is, it may   have multiple edges,
and also loop edges.} 
is a tree. See Figure~\ref{Wtree} for an example of 
a chord diagram $D$ whose associated graph is the tree \Tree~.  

\begin{figure}[h]
\begin{center}
$W_{2}(\ACE)=W_{0}(\ACEi)=1$
\caption{\label{Wtree}A degree $2$ chord diagram $D$ with $W_{2}(D)=1$.   }
\end{center}
\end{figure}
Thus, we come to the following result.
\begin{lemma}\label{2.3} For chord diagrams on $m$ circles, the
  Alexander-Conway weight system satisfies $W_{i}=0$ for
   $i\leq m-2$. Moreover,  $W_{{m-1}}$ takes the value $1$ on
   precisely those chord 
  diagrams whose associated graph is a spanning tree on the complete
   graph  $K_m$, and $W_{{m-1}}$ is zero on all other chord diagrams.
\end{lemma}

On the other hand,
the linking number $\ell_{ij}$ is a finite type invariant of order $1$ 
whose weight system 
is the linear form dual to the chord diagram having just
one chord connecting the $i$th and $j$th circle. It follows that
the Kirchhoff polynomial $\D_m$ (see~(\ref{eq:kirch_pol})) in the
linking numbers $\ell_{ij}$ is a finite type invariant of order $m-1$ 
and that its weight system, by our lemma, is equal to $W_{{m-1}}$. Thus, we
have proved Theorem~\ref{Hoste.thm}  on the level of weight systems.  

\begin{remark}
\rm This actually implies the
theorem, since the Alexan\-der-Conway polynomial is (almost) a 
 canonical invariant,
\emph{i.e.}\ it can be recovered from its weight system by the
 Kontsevich integral. 
We will discuss this 
in a more general situation
 in the proof of Proposition~\ref{5.3}. 

This proof of Theorem~\ref{Hoste.thm}
can be given without 
mentioning the theory of finite type invariants.
Indeed, 
it can be reformulated as a proof
by induction  on the number of 
components using the skein relation~(\ref{skein}) (see the
original arguments in~\cite{Ha} and~\cite{Ho}).  
The point of our formulation is that it 
can be generalized
from ordinary linking numbers to higher order Milnor
numbers, where finite type invariants technology will 
play a crucial role.
\end{remark}

\section{A Vanishing Lemma}\label{sec.main}

In this section, we generalize Lemma~\ref{2.3} 
and show that the 
Alexander-Conway weight system vanishes also on some classes of diagrams 
of degree higher than $m-2$.

We now allow not just chord diagrams, but also diagrams with internal
(trivalent) vertices. Let us briefly review the relevant concepts and fix our
notation. (For more on diagrams and finite-type invariants see 
\emph{e.g.}~\cite{BN}.) 

A {\em diagram} $D$ on an oriented $1$-manifold $X$ is a finite
uni-trivalent graph $\Gamma$ 
whose univalent vertices are attached to (the interior of) $X$, and so
that  each  trivalent vertex  is equipped with a cyclic  
ordering 
of the three half-edges meeting at that vertex.
It is assumed that every component of $\Gamma$ has at
least one univalent vertex.   

The \emph{degree} of a diagram $D$ is one half of the total number of
vertices of $\Gamma$ (both univalent and trivalent vertices count). We
denote by $\A_d(X)$ the $\Q$-vector space 
spanned by degree $d$ diagrams on $X$, modulo IHX, AS, and STU
relations (see Figure~\ref{IHX}).  
 
\begin{figure}[h]
\begin{center}
IHX:\ \  \IHX\\
AS:\ \ \ASi\\ 
STU:\ \ \STU
\end{center}
\caption{\label{IHX} IHX, AS, and STU relations. Orientations at $3$-valent 
vertices are those induced from the orientation of the plane.
Four-valent vertices are an artifact of the planar pictorial representation
and should be ignored.}
\end{figure}

Traditionally, the components of $X$  are 
referred to as {\em solid}, while the components of $\Gamma$ are
referred to as {\em   dashed}, since the edges of $\Gamma$ are usually
represented pictorially  by  dashed lines. Thus, for example, elements
of the space $\A_d(\amalg_m S^1)$ are referred to as
 (linear combinations of) 
diagrams of degree $d$ on $m$ solid circles. 

Let us agree that by the {\em   components} of a diagram $D$, 
we always mean the components of its {\em  dashed} part $\Gamma$. 
A component is called a {\em tree component} if its 
underlying graph is a tree. Note that the degree of a tree component
is one less than the number of its univalent vertices.  
For example, a chord component \Fi has degree one, 
and a $\Y$-shaped component \MILN has degree two.

\medskip

Recall that the weight system of the coefficient $c_n$ of the 
Alexander-Conway polynomial is denoted by $W_n$. This  
is a linear form on the $\Q$-vector space $\A_n(\amalg_m S^1)$. 
We will simply write $W$ for $W_n$, if the degree is clear from the
context.  

Note that~(\ref{eq:parvanish}) gives
\begin{equation} \label{eq:parity}
W_n=0
\ \ \text{if \ $n\equiv m$ {\em mod} $2$}
\end{equation}
(this also follows immediately from~(\ref{smooth})
and~(\ref{compW})).

\begin{proposition}[{\bf Vanishing Lemma}] \label{2.1} 
Let $D$ be a  degree-$d$ diagram 
 on \hbox{$m\geq 2$} solid circles, such that $D$ has no tree
 components of degree $\leq n-1$. 
If $d\leq n(m-1)+1$,  then $W(D)=0$ unless $D$ has exactly $m-1$
 components, each of which is a tree of degree $\geq n$.
\end{proposition}
\begin{corollary}\label{r4.2}   
If a diagram $D$ satisfies the  hypotheses of the Vanishing Lemma, 
and 
 $W(D)\neq 0$, then either 
 $d=n(m-1)$ or $d=n(m-1)+1$. 
Moreover,  the only two possibilities are  the following:   
\begin{itemize}
\item If at least one of the numbers $n$ and $m$ is odd, then
$d=n(m-1)$ and each of the $m-1$ components of $D$ is a tree 
of degree $n$. 
\item If both $n$ and $m$ are even, then  $d=n(m-1)+1$ and $D$ has 
$m-1$ components of which $m-2$ are trees of degree $n$ 
   and one is a tree of degree $n+1$. 
\end{itemize}
  \end{corollary}
\noindent
\emph{Proof of  Corollary~\ref{r4.2}.} This follows from 
Proposition~\ref{2.1} and~(\ref{eq:parity}),
since 
$D$ must be of degree $d\ge n(m-1)$.
\eproof 
\medskip

In the proof of the Vanishing Lemma~\ref{2.1}, it will be convenient
to use the following 
terminology.
If the Alexander-Conway weight system takes equal values  
$$ W(D)=W(D')$$
on two elements  
$$
D\in \A_d(X) \text{\ and \ } D'\in \A_{d'}(X')
$$ 
(possibly of different degrees, or maybe even with $X\ne X'$, in which 
case the correspondence between $X$ and $X'$ should be clear from the
context), then we call $D$ and $D'$  \emph{equivalent} and write this
as   
$$ D \equiv D'~. $$

We begin with general formulas that allow to evaluate the
Alexander-Conway weight system on $\A_d(\amalg_m S^1)$ for all $d$ and
$m$. 

\begin{lemma}[\cite{FKV}]\label{ACRel}
The Alexander-Conway weight system satisfies the following 
       reduction relations:

\begin{equation}\label{Relb}
\Relb \ \ \ \equiv \ \ \Relbi
\end{equation}
\begin{equation} \label{Rela}
\Rela \ \ \ \equiv \ \ -2 \ \ \ \Relai
\end{equation}
\begin{equation}\label{Relc}
\Relc \ \ \ \equiv \ \ 0
\end{equation}
\begin{equation}\label{Reld}
\Reld \ \ \ \equiv \ \ 0
\end{equation}
\begin{equation}\label{Rele}
\Rele \ \ \ \equiv \ \ - \frac{1}{2} \Bigl( \ \ \Relei +\ \put(0,-.4)
{\rotatedown\Relei} \ \ \ \ \ \ \ \ -\ \Releii -\ \Releiii \ \ \Bigr) 
\end{equation}
\begin{equation} \label{Relconseq}
\Stepdi \ \ \equiv \Stepdii \ \ + \Stepdiia \ \ -\Stepdiib\ \ - \Stepdiic
\end{equation}
\end{lemma}

\proof
Relations~(\ref{Rela})---(\ref{Rele}) have been proved in~\cite{FKV}. 
They follow from the  connection between 
the Alexander-Conway weight system  $W$ and the Lie superalgebra
$gl(1|1)$.  Equation~(\ref{Relb}) is equivalent to the fact
that the invariant metric on $gl(1|1)$ needed to produce $W$ 
is given by the trace in the defining representation.
Finally, equation~(\ref{Relconseq}) is a corollary
of~(\ref{Rele}) and~(\ref{Rela}). 
\eproof

\begin{remark}
\rm In this paper, we will mainly need to apply these relations to
    diagrams where all univalent vertices lie on 
    solid  circles. For such  diagrams, these relations can also
    be deduced from the smoothing relations~(\ref{smooth})
    and~(\ref{compW}), together with the STU relation.
\end{remark}

\begin{lemma} \label{2.2} Let $D$ be a diagram on $m\geq 2$ solid
circles with $a+b$
components that consist of
\ $a$ \ chords \Fi and \ $b$ \
   $\Y$-shaped components \MILN .
If $a+b<m-1$, then $W(D)=0$.  
\end{lemma}
(Note that if $b=0$, this  
has already been proved in Lemma \ref{2.3}.)

\begin{figure}[h]
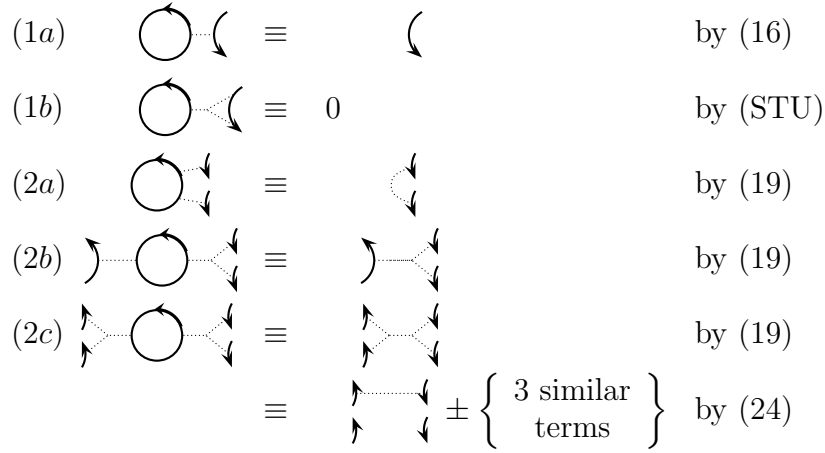

\begin{center}
$$\begin{array}{lclll}
(1a)&\Stepa  &\equiv& \ \ \ \ \ \ \ \   \Stepai \  \ &\text{by
  (\ref{smooth})} \\
(1b)&\Stepaa   &\equiv& \ 0  \ &\text{by (STU)} \\
(2a)&\ \ \Stepb  &\equiv& \  \Stepbi \ \ \ \  & \text{by (\ref{Relb})} \\
(2b)&\ \ \ \Stepc  &\equiv &\  \Stepci \ \ \ \  &\text{by (\ref{Relb})} \\
(2c)&\ \ \Stepd  &\equiv &\  \Stepdi \ \ \  &\text{by (\ref{Relb})} \\
&&\equiv& \Stepdii \ \pm \left\{ 
\begin{array}{c} \text{3 similar}\\
\text{terms}
\end{array}\right\}&\text{by (\ref{Relconseq})}\\ 
\end{array}$$
\caption{\label{onevert} Reducing  $D$ to $D'$. There are five
  cases to consider.  Note that in the case (2c), $ a$ increases by
  $1$, but $b$ decreases by $2$. In all cases, both $m$ and $a+b$
  decrease by $1$. } 
\end{center}
\end{figure}

\proof The proof is by induction on $m$. If $m=2$, then $a=b=0$  and
so $W(D)=0$ by (\ref{unlink}). Now assume $m>2$. Since $a+b<m-1$, the
number of univalent vertices of $D$, which  is $2a+3b$, is strictly
less than $3m$. Therefore at least one of the solid circles has at
most $2$ univalent vertices on it. If this solid circle has 
no univalent vertices, then $W(D)=0$ by (\ref{unlink}). If it has one 
or two univalent vertices, we apply 
the reduction relations of Lemma~\ref{ACRel}  to get $D\equiv D'$, 
where $D'$ is a diagram on $m-1$ circles with $a'$ chords and $b'$
$\Y$'s, and $a'+b'=a+b-1$ (see Figure \ref{onevert}). Thus
$W(D)=W(D')=0$, where the second equality is  by the induction
hypothesis. This proves the lemma. 
\eproof
 
\begin{lemma} \label{lemred} A tree component of odd (resp., even)
degree is equivalent to a linear combination of chords 
(resp., $\Y$'s).  
\end{lemma} 

\begin{figure}[h]
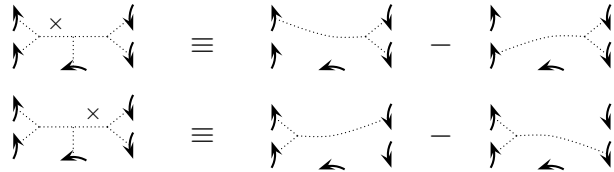

\begin{center}
$$\begin{array}{ccc}
\Treeqa &\equiv& \Treeqi -\Treeqii\\
\Treeqb       &\equiv& \Treeqiii -\Treeqiv
\end{array}$$
\caption{\label{Tree4} Two ways to reduce a tree of degree $4$. We
  apply (\ref{Rele}) at the marked edge, and then apply (\ref{Rela})
  and (\ref{Relc}).}
\end{center}
\end{figure}

\begin{figure}[ht]
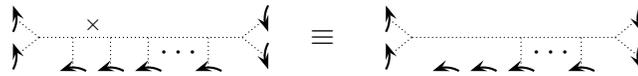

  \begin{center}
    \treea \ \ \ $\equiv$ \treeb
  \caption{\label{reducetrees}How to reduce trees of degree $\geq
    5$. Because of (\ref{Relc}), only one term survives after 
    applying (\ref{Rele}) at the marked edge. }
\end{center}
\end{figure}

\proof Equation~(\ref{Relconseq}) in 
Lemma~\ref{ACRel} shows that a tree of degree $3$ is equivalent
to a linear combination of chord diagrams. Figure~\ref{Tree4} shows
that a tree of degree $4$ is equivalent 
to a linear combination of $\Y$'s. Trees of higher degree are either 
equivalent to zero by (\ref{Reld}) or can be reduced recursively to
trees of degree $3$ or $4$ by the move shown  in
Figure~\ref{reducetrees}. \eproof 
 
\noindent
\emph{Proof of  Proposition~\ref{2.1}.}
Let 
$d\leq n(m-1) +1$ be minimal such that there
exists a diagram $D\in \A_d(\amalg_m S^1)$ without tree components of
degree $\leq n-1$ and $W(D)\neq 0$. 

\begin{figure}[ht]
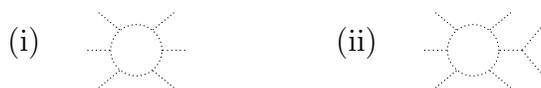

 \begin{center}
   (i) \Wheel  \ \ \ \ \ \ \ \ \ \ \ (ii) \Wheelb
   \caption{\label{wheel}(i) A wheel with $6$ legs, and (ii) a component 
  which is neither a   wheel nor a tree.}
 \end{center}
\end{figure}

First, we will show that every component of $D$ is a tree. Indeed,  
assume that a component  of $D$ is a wheel (see Figure~\ref{wheel}(i))     
with $k\geq 2$ legs. (A  wheel with only one leg is zero by the STU
relation.) If $k=2$,  then we can remove the wheel
using~(\ref{Rela}). The diagram $D'$ thus obtained has degree $d-2$,
and still has no trees of degree $\leq n-1$. 
But $W(D')=-2 W(D)\neq 0$, which contradicts the minimality 
of $d$. Similarly, if the wheel  
component has $k\geq 3$ legs, then we can apply~(\ref{Rele})
and~(\ref{Rela}) as in Figure~\ref{reducetrees} and get a wheel 
with $k-2$ legs.  Again, the diagram $D'$ thus obtained has smaller
degree  and hence contradicts the minimality of $d$.  Therefore no
component of $D$ can be a wheel. 
Furthermore, if a component of $D$ is
neither a wheel nor a tree
(see \emph{e.g.}\ Figure~\ref{wheel}(ii)), then it has 
a trivalent vertex all of whose neighboring vertices are also 
trivalent vertices, hence $W(D)=0$ by~(\ref{Reld}).

Hence every component of $D$ is a tree as claimed.

Now let \ $a$ \ (resp.\ $b$) be the number of trees of odd (resp.\ even)
degree in $D$. By hypothesis, there are no trees of degree $\leq n-1$, 
hence $d\geq n(a+b)$. Since  $d\leq n(m-1)+1$, it follows that  
$a+b\leq m-1 + \frac{1}{n}$. 
Since the case $n=1$ has already been dealt with in 
Lemma~\ref{2.3}, we may assume $n\geq 2$.
This gives
$a+b\leq m-1$.
  
Reducing the tree components to chords or $\Y$'s by 
Lemma~\ref{lemred}, we see that our diagram  $D$ is equivalent to a
linear combination of diagrams each of which has exactly $a$ chords
and $b$ $\Y$'s. If $a+b<m-1$, then $W(D)=0$ by Lemma~\ref{2.2}. Thus
we must have $a+b = m-1$, and $D$ has exactly $m-1$ components   each of 
which is a tree of degree $\geq n$.

Thus we have proved the result for the diagram $D$. 

We have assumed that $D$ had minimal degree, but since this degree is
already $\geq n(m-1)$, and $W_d\neq 0$ 
implies $W_{d+1}=0$ by~(\ref{eq:parity}), 
diagrams of higher degree  will either be killed by $W$ or do not
satisfy the hypotheses any more.  

This completes the proof.
\eproof

\section{%
The first non-vanishing Milnor invariants}
\label{sec:milnor}

In \cite{M1,M2}, Milnor introduced integer-valued 
invariants $\mu_{i_1 i_2 \ldots i_k}(L)$ of 
oriented links, known as Milnor higher linking numbers 
(the number $k$ of indices is called the {\em length} 
of a Milnor invariant).

In fact, these invariants are not 
universally defined, 
that is,
unless the lower order (\emph{i.e.}\ of lower length) invariants  
vanish, they are either not defined or have indeterminacies. 

The first Milnor invariants are  the ordinary linking numbers 
$\mu_{ij}(L)=\ell_{ij}(L)$   which are always defined.  If they all
vanish, then the triple linking 
numbers $\mu_{ijk}(L)$ are well-defined. In general, if all invariants      
of length $\leq n-1$ are defined and are zero, then all invariants of 
length $n$ are defined. 

Habegger and Lin~\cite{HL1} introduced the 
philosophy that Milnor's invariants are actually invariants of string
links which are always defined. 
The first non-vanishing Milnor
invariants of a link $L$ are well-defined and equal to the lowest 
degree   non-vanishing invariants of any string link representative of
$L$, and the indeterminacies of invariants  of higher degree
come precisely from the indeterminacy 
of representing a link as the closure of a string link.  From this
point of view, 
Bar-Natan~\cite{BN1} and Lin~\cite{L} showed that Milnor's invariants 
of string links are of finite type. 

Here, the \emph{degree} of an invariant is one less than 
its length (the number of its indices).
For example, linking numbers have degree one, and
triple linking numbers have degree two.

The universal finite type invariant of a string link $\s$ is its Kontsevich
integral $Z(\s)$. It was shown in \cite{HM} how to compute the Milnor
invariants of $\s$ from the tree part of $Z(\s)$. This leads to the following
description of the first non-vanishing invariants.  

By a {\em labelled diagram}\footnote{These diagrams are 
sometimes also called Chinese characters,  Feynman diagrams,  web diagrams,  
Jacobi diagrams, 
\emph{etc.}}  
of degree $n$ we mean a uni-trivalent graph with $2n$ vertices
with a cyclic ordering of the three half-edges meeting at each
trivalent vertex 
and whose univalent vertices are labelled by elements 
     of the set $\{1,\ldots, m\}$. 
 We denote by $C_n(m)$ the $\Q$-vector space of labelled diagrams 
of degree $n$
modulo AS and IHX relations (see Figure~\ref{IHX} in
Section~\ref{sec.main}).

We will be interested in the subspace 
$$C^t_{n}(m) \subset C_n(m)$$
generated by those labelled diagrams whose
underlying graphs  are trees. 
In other words,
$C^t_{n}(m)$ is just the space of  
tree diagrams
of degree $n$ with labels from
$\{1,\ldots, m\}$ on their univalent vertices, modulo  AS and IHX 
relations.   
 
Denote by $\L_n(m)$ the set of 
 oriented (framed) links $L$ with $m$ (numbered) components, such that all
 Milnor invariants of $L$ of degree less than $n$
vanish. This defines a decreasing filtration 
$$ \L_1(m) \supset \L_2(m) \supset \ldots \L_n(m) \supset \ldots $$
on the set of
 links. (Note that every link lies in  $\L_1(m)$.) For links in
 $\L_n(m)$, Milnor invariants of degree $n$ are well-defined.   

\begin{remark} \rm 
We 
consider unframed links as elements of $\L_1(m)$ by equipping them
with the zero framing. 
However, for our purposes the choice of framing does not really  
affect anything, since the only framing-dependent Milnor invariants
are the self-linking numbers $\ell_{ii}(L)$ and the Alexander-Conway 
polynomial ignores framing. 
\end{remark}

As explained in~\cite{HM}, the set of all degree-$n$ 
Milnor invariants of a link $L\in \L_n(m)$ 
can be encoded by a universal invariant  
\begin{equation} \label{eq:xi}
        \xi_n : \L_n(m) \to      C^t_{n}(m)
\end{equation}
with values in the space of tree diagrams of degree $n$.

\begin{definition}\rm
We call $\xi_n(L)\in C^t_{n}(m)$ the \emph{universal degree-$n$ Milnor
invariant} of a link  $L\in \L_n(m)$.
\end{definition}

The cases of $n=1$ and $n=2$
can be described very explicitly, as follows. 
Note that $$C^t_1(m)\cong S^2V$$ and  $$C^t_2(m)\cong \Lambda^3V,$$
where $V$ is a vector space of dimension $m$ 
with a fixed basis corresponding to the labels $1,2,\ldots,m$.
Thus $C^t_1(m)$ has a basis represented by `struts'  
\vspace{4pt}

\centerline{${\scriptstyle i}\ \Fi\ {\scriptstyle j}$}

\vspace{4pt}

The coefficient of this strut in $\xi_1(L)$ is equal to
the linking number $\ell_{ij}(L)$ between the $i$th and $j$th
component of $L$.  (We assume $i\neq j$ here.) 
We can therefore view $\ell_{ij}$ as the linear form on 
$C^t_1(m)\cong S^2V$ which sends $\xi_1(L)$ to $\ell_{ij}(L)$.   

If $L\in \L_2(m)$ (that is, if all linking numbers of $L$ vanish),
then $\xi_1(L)=0$ and $\xi_2(L)$ is well-defined. It lies in
$C^t_2(m)\cong \Lambda^3V$ which has a basis 
represented by (antisymmetric) $Y$-shaped labelled diagrams 

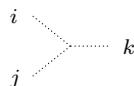
\centerline{
\MILNijk
} 
\noindent
The coefficient of this diagram in $\xi_2(L)$ is the Milnor triple
linking number $\mu_{ijk}(L)$. In other words, $\mu_{ijk}$ is the
linear form 
on $C^t_2(m)\cong \Lambda^3V$ which sends $\xi_2(L)$ to $\mu_{ijk}(L)$. 

Similarly, for links $L\in \L_n(m)$, we have $\xi_i(L)=0$ for $i<n$,
and $\xi_n(L)$ is well-defined. 
Milnor invariants  of degree $n$
can be considered as linear forms on $C^t_{n}(m)$ sending 
$\xi_n(L)$ to the numerical invariants, \emph{e.g.}\  to 
$\mu_{i_1 i_2 \ldots i_{n+1}}(L)$. 
In particular, the dimension of $C^t_{n}(m)$ is equal to the number of 
linearly independent Milnor invariants of degree $n$ for links in
$\L_n(m)$.\footnote{This number, which was first computed by
Orr~\cite[Theorem 15]{O}, is equal to $mN_n(m)-N_{n+1}(m)$, where 
$\sum_i  N_i(m)t^i$ is the Hilbert series of the free Lie algebra on
$m$ letters. See also~\cite[Section 8]{HM}.}                 
We refer the reader to~\cite[Sections 5 and 6]{HM} 
  for an explicit description 
of the linear forms 
giving
numerical Milnor invariants.
(We will not need it in this paper.)
The general case is
slightly more complicated than 
the cases of degrees one and two,
because the space $C_n^t(m)$ for $n\ge 3$ does not have a canonical basis.
\medskip

One motivation for encoding  the first non-vanishing Milnor invariants
of a link $L\in \L_n(m)$ into the universal invariant $\xi_n$
is that $\xi_n(L)$ can be identified with the degree-$n$ term $Z_n^t(\s)$ in 
the tree part of the Kontsevich integral of a  string link $\s$ whose
closure is $L$.  Since this identification will be used
in the next section, we will now review it here (see~\cite{HM}
for details).

The Kontsevich integral $Z(L)$ of an $m$-component (framed) link $L$
lies in the space $\A(\amalg_m S^1)$ of diagrams on $m$ solid
circles. Thus 
it
is a power series $$Z(L)=\sum_{n\ge 0} Z_n(L),$$ 
whose $n$th term $Z_n(L) \in \A_n(\amalg_m S^1)$ lies in   
the space of diagrams of degree $n$.

Suppose now that 
$L$ is the closure of some string link $\s$. 
Its
Kontsevich integral  
is a power series $$ Z(\s)=\sum Z_n(\s),$$
 whose $n$th term
$Z_n(\s)$ lies in 
the space 
$$\A_n(m):=\A_n(\amalg_m \I)$$ 
of degree
$n$ diagrams on $m$ solid intervals. 

Denote by $\A^t(m)$ the quotient space of $\A(m)$
modulo
diagrams containing non-simply connected (dashed) components. 
(The superscript $t$ stands for tree diagrams, which generate
$\A^t(m)$ multiplicatively.) 
By definition, the \emph{tree part} of the Kontsevich integral 
$$Z^t(\s)=\sum Z^t_n(\s)$$
is the image of  $Z(\s)$ in $\A^t(m)$.

For every string link $\s$
we have 
        $$Z_0^t(\s)=1,$$
where $1$ is the diagram consisting of the trivial string link without any  
dashed part. 
It is shown in \cite{HM} that if $L\in \L_n(m) - \L_{n+1}(m)$ 
({\em i.e.}, 
the first non-vanishing Milnor invariants of $L$ occur precisely in degree
$n$), then the first non-vanishing term of $Z^t(\s)-1$ also occurs in degree
$n$. (In other words, the {\em $Z^t$-filtration degree} of $\s$ is equal to
the {\em Milnor filtration degree} of $L$.) Moreover, since $Z^t(\s)$ is
group-like, $Z^t_n(\s)$ lies in the primitive part of $\A_n^t(m)$, which is
naturally identified with $C_n^t(m)$. Under this identification,
$Z^t_n(\s)$ becomes $\xi_n(L)$, 
the universal degree-$n$ Milnor invariant of $L$. In particular,  
all first non-vanishing Milnor invariants of $L$ can be computed
from $Z^t_n(\s)$ by an explicit formula given in~\cite[Theorem~6.1]{HM}. 

The identification of the primitive part of $\A_n^t(m)$ with
$C^t_{n}(m)$ is  very straightforward. As in~\cite[Section 4]{HM}, let
us denote the primitive 
part by $\P_n^t(m)\subset \A_n^t(m)$. It is the subspace generated by
diagrams with connected dashed parts. The dashed part in such a diagram is a
tree, hence after removing the $m$ solid intervals and retaining only their
labels in $\{1,\ldots,m\}$ on the univalent vertices, we get an element of
$C^t_{n}(m)$. 
The relations in $\A^t(m)$ together with the STU relation imply that
this gives a well-defined 
map $$\P_n^t(m) \mapright\sim C^t_{n}(m),$$
and moreover, this map is an isomorphism.

\section{A formula for the coefficient $c_{n(m-1)}$ }
\label{sec:coeff}

In this section, we show that if the first non-vanishing Milnor
invariants of an $m$-component link $L$ have degree $n$ 
({\em i.e.}\ if $L\in \L_n(m)$), then 
the coefficients of its Alexander-Conway polynomial $\nabla_L(z)$  
vanish in degrees smaller than $ n(m-1)$ and
$$ c_{n(m-1)}(L) = F_m^{(n)}(\xi_n(L)),$$
 where $\xi_n(L)\in C^t_n(m)$ 
is the universal degree-$n$ Milnor invariant~(\ref{eq:xi}) 
of  $L$, and $F_m^{(n)}$ is a homogeneous polynomial of degree $m-1$
on  $C^t_n(m)$. 
Using results and techniques of~\cite{HM}  this 
will follow from  the Vanishing Lemma~\ref{2.1}.      

The homogeneous polynomial $F_m^{(n)}$ will be given in terms of the
corresponding multilinear form 
$\widetilde  F_m^{(n)}$ 
on $C^t_n(m)$  defined as follows. 

We call a {\em lift} of a labelled diagram 
$\gamma\in C_d(m)$ any 
diagram
              $$D\in     \A_d(\amalg_m S^1)$$
on $m$ solid circles
 whose dashed components are the same as the components of $\gamma$ and
for each $i$, the univalent vertices 
of $\gamma$
 with label $i$  are  attached  in $D$ 
to the $i$th solid circle.  

 Given  tree diagrams $$\xi^{(1)}, \ldots, \xi^{(m-1)}\in C^t_n(m),$$
 consider their  disjoint union  $\gamma$ 
as an element of $C^t_{n(m-1)}(m)$, and let  $D$  be
 a lift of 
$\gamma$ to $\A_{n(m-1)}(\amalg_m S^1)$. 
Then    we set
\begin{equation}\label{eq:fmtilde}
\widetilde F_m^{(n)}(\xi^{(1)}, \ldots,\xi^{(m-1)})=W_{n(m-1)}(D),
\end{equation}
where $W_{n(m-1)}$ is the Alexander-Conway weight system. 

\begin{lemma}\label{Fdef} 
The multilinear form $\widetilde F_m^{(n)}$ is well-defined. 
\end{lemma}
\proof
If $D'$ is another lift of $\gamma$,  
then $D'$ is obtained from $D$ by permuting the
cyclic order of the univalent vertices on every 
component of $\amalg_m S^1$. 
All such permutations can be obtained by composing simple transpositions
of two adjacent vertices on a single circle. Therefore, by 
applying the STU relation, we see 
that $D$ and $D'$ differ by diagrams containing non-simply connected  
components  or tree components  of degree $\ge n+1$. Since
$W_{n(m-1)}$ is zero on such diagrams by the 
Vanishing Lemma~\ref{2.1}, the result follows. 
\eproof

Now, for a tree diagram $\xi\in C^t_n(m)$  we define 
\begin{equation}\label{eq:fn}
        F_m^{(n)}(\xi):= \frac{1}{(m-1)!}
                         \widetilde F_m^{(n)}(\xi,\ldots,\xi),
\end{equation}
so that $\widetilde F_m^{(n)}$ is the polarization of $F_m^{(n)}$.

\begin{remark}\label{r6.2}
\rm
Note that $F_m^{(n)}=0$ if both $m$ and $n$ are even. This
follows immediately from~(\ref{eq:parity}).
\end{remark}

\begin{proposition}\label{5.2} 
Let $L\in \L_n(m)$, 
\emph{i.e.}\ all Milnor invariants of $L$ of degree $\le n-1$ vanish. 
\begin{itemize} 
\item[(i)] The Kontsevich integral $Z(L)$ can be written as a linear
  combination of diagrams none of which contains a tree component of degree
  $\leq n-1$.  
 \item[(ii)] Modulo diagrams containing non-simply connected
  components or tree components of degree $\ge n+1$, we have
  $$Z_{n(m-1)}(L))={1\over (m-1)!}\, \xi_n(L)^{m-1},$$  
where $\xi(L) \in  C^t_n(m)$ is the universal degree-$n$ Milnor
invariant~(\ref{eq:xi}) of $L$. 
\end{itemize}
\end{proposition}

Here, by abuse of notation, $\xi_n(L)^{m-1}$ is considered as 
an element of $\A_{n(m-1)}(\amalg_mS^1)$. 
This is acceptable, since any two lifts of
$\xi_n(L)^{m-1}$ to the solid circles differ by diagrams containing
non-simply connected components or tree components of degree 
$\ge n+1$
(see the proof of Lemma~\ref{Fdef}). 

\proof
Part (i) is Corollary~9.3 of~\cite{HM}. The proof of~(ii)
is based on a similar argument.   
Here, we only sketch the idea and refer to \cite{HM} for more
details. One has 
$$
\log \,Z^t(\sigma)
=Z^t_n(\sigma) +O(n+1),
$$ 
where $O(n+1)$ stands for terms of degree $\geq n+1$. 
It follows that modulo trees of degree $\geq n+1$
we have
$$
Z_{n(m-1)}^t(\s)={1\over (m-1)!}\, (Z^t_n(\sigma))^{m-1}~.
$$ 
Now
$Z^t(L)$ is obtained from $Z^t(\sigma)$ by a procedure involving 
a special element $\nu$ in the algebra of diagrams
(it is related to  the Kontsevich integral of the unknot). 
Since $\nu$ does not contain trees, this procedure does not introduce any new
trees. Since $Z^t_n(\sigma)$ is identified with $\xi_n(L)$, as explained in
the previous section, this gives~(ii).  \eproof

\begin{proposition} \label{5.3} 
Let $L\in \L_n(m)$, and let 
$\xi_n(L)\in   C^t_n(m)$   
be its universal degree-$n$ Milnor invariant.
Then  for the coefficients 
$c_i(L)$ of the Alexander-Conway
polynomial 
$\nabla_L(z)=\sum_{i\ge 0}c_i(L)z^i$ 
we have 
\begin{itemize} 
        \item[(i)] \ $c_i(L)=0$ for $i<n(m-1),$
\item[(ii)] \ $c_{n(m-1)}(L)= F_m^{(n)}(\xi_n(L)).$
\end{itemize}
\end{proposition}

\proof Consider the following renormalized version of the
Alexander-Conway polynomial $\nabla_L(z)$:  
\begin{equation}
  \label{eq:tn}
 \widetilde \nabla_L(z) =
 \frac{z}{e^{z/2}-e^{-z/2}} \nabla_L(e^{z/2}-e^{-z/2})
=\sum_{n\geq 0} \tilde c_n(L) z^n~.
\end{equation} 
Note  that the coefficient $\tilde c_n$ is a  finite type invariant of
 order  $n$. 
Its weight system is equal to $W_n$ (the weight system of $ c_n$).   

It was shown in \cite{BNG} that  $\tilde c_n$ is a  \emph{canonical
invariant}, \emph{i.e.}\
it can be recovered from its weight system by the Kontsevich integral:
\begin{equation}
  \label{eq:2}
\tilde c_n(L) = W_n(Z_n(L))~.
\end{equation}

  Now assume that  $L\in \L_n(m)$, \emph{i.e.}\ the first non-vanishing
Milnor invariants of $L$ have degree $n$.
By the Vanishing Lemma~\ref{2.1}  and  Proposition~\ref{5.2}, 
  non-vanishing coefficients of $\widetilde \nabla_L(z)$   
can only occur in degrees   $\ge n(m-1)$, and 
\begin{align*}
\ \tilde c_{n(m-1)}(L)=& W_{n(m-1)}(Z_{n(m-1)}(L))=  {1\over (m-1)!}\,
W_{n(m-1)}(\xi_n(L)^{m-1})\\ 
=&\ {1\over (m-1)!}\widetilde F^{(n)}_m(\xi_n(L),\ldots, \xi_n(L)) =
F^{(n)}_m(\xi_n(L))~. 
\end{align*}

Since the first non-vanishing coefficient of $\nabla_L(z)$ is equal to
the one of $\widetilde \nabla_L(z)$, 
this gives the desired result.
\eproof 

\section{Expressing $c_{n(m-1)}$ via spanning-tree polynomials}
\label{sec:fm}

In this section we will rewrite the polynomials $F_m^{(n)}$ introduced
in the previous section in terms of the spanning-tree polynomials 
$\D_m$ and  $\P_m$ (see equations~(\ref{eq:kirch_pol}) and~(\ref{eq:pm})). 

This will furnish the proofs of the tree-sum formulas of
Theorems~\ref{thm:alex-trees} 
and~\ref{thm:general} for the 
coefficient $c_{n(m-1)}$ of the Alexander-Conway polynomial.

\subsection{The case $n=1$.} 

\

The case $n=1$ is given by  Theorem~\ref{Hoste.thm}. With the
identification  $C_1^t(m)= S^2V$, the polynomial $F_m^{(1)}\in
S^{m-1}(C_1^t(m))^*$ is just the Kirchhoff polynomial $\D_m$ 
in the linking numbers $\ell_{ij}$ (see~(\ref{eq:kirch_pol})).
 This is also clear from the discussion  in Section~\ref{sec:ac}
(see Lemma~\ref{2.3}).

\subsection{The case $n=2$.}

\

As already observed in Remark~\ref{r6.2}, we have $F_m^{(2)}=0$ 
if $m$ is even. 

Let us now assume that $m$ is odd. 
Since $F_m^{(2)}$ is a homogeneous polynomial 
of degree $m-1$ on $C^t_2(m)\cong \Lambda^3 V$,
it can be written as a sum of monomials of degree $m-1$
in the variables $\mu_{ijk}$ (which represent the standard basis  of
the dual space  $C^t_2(m)^*$). 
The coefficient in $F_m^{(2)}$ 
of each  monomial in the $\mu_{ijk}$'s
can be computed by using the Alexander-Conway weight system. 

This gives the following characterization of $F_m^{(2)}$. 
We will write
$$\mu_{ijk}=v_i\wedge v_j\wedge v_k~,$$ 
where 
$v_1,\ldots, v_m$ is the standard basis of $V^*$. 
This allows us to consider $F_m^{(2)}$ as an 
expression in the variables $v_i$: 
 $$ F_m^{(2)}=F_m^{(2)}(v_1,v_2,\ldots, v_m)~.$$ 
Note that we know  already that this expression  
is invariant under any permutation of the $v_i$'s. 
(This follows from  Proposition~\ref{5.3}, since the
coefficient $c_{n(m-1)}(L)$ does not depend on an 
ordering of the components of $L$.) 

\begin{proposition} \label{recursion}  
Assume   that $m\geq 3$ is odd. Then the polynomial
  $F_m^{(2)}$ is characterized by the following properties. 
(For   simplicity, we write  $F_m$ 
instead of
$F_m^{(2)}$.)
\begin{itemize}

\item[(i)] $F_3=\mu_{123}^2$

\item[(ii)] If a monomial 
$\prod_{\alpha} \mu_{i_\alpha j_\alpha k_\alpha}$ 
  occurs with non-zero coefficient in $F_m$, then there exists
  $p\in\{1,2,\ldots, m\}$ such that $p$ occurs exactly twice in the
list of indices $i_1,j_1,k_1,i_2,\ldots j_{m-1},k_{m-1}$. 

\item[(iii)] $F_m$ satisfies the following relations
(and all the relations obtained from them by permuting the indices
$1,2,\ldots,m$):
\begin{align*}
\left[\frac {\partial^2 \, F_m}{\partial \mu_{123}^2}\right]_{v_1=0}
&= 2 F_{m-2}( v_2+v_3,\ldots),\\
\left[\frac {\partial^2 \, F_m}{\partial \mu_{123}\, \partial
    \mu_{124}}\right]_{v_1=0}&= F_{m-2}( v_2+ v_3, v_4, \ldots) +
F_{m-2}(v_2+v_4,v_3, \ldots)\\ 
& \ - F_{m-2}( v_3+v_4,v_2, \ldots),\\
\left[\frac {\partial^2 \, F_m} {\partial\mu_{123}\,
\partial\mu_{145}}\right]_{v_1=0} &= F_{m-2}( v_3+v_4,v_2, v_5, \ldots)
+ F_{m-2}(v_2 +v_5, v_3, v_4, \ldots)\\ 
& \ -F_{m-2}( v_2+v_4,v_3, v_5, \ldots) - 
 F_{m-2}( v_3 +v_5, v_2, v_4, \ldots) ~.
\end{align*} 
(Here
the dots stand for the $v_i$ with indices  which do not appear in 
the left hand side; for example,  in the first equation, 
the dots
mean
$v_4,v_5, \ldots, v_m$).
 \end{itemize}
\end{proposition}

\proof The key point is to observe that
monomials in the $\mu_{ijk}$'s 
are dual to diagrams consisting of $\Y$'s
glued to the $m$ solid circles. 
With this in mind, statement~(ii) is a 
reformulation of the following
fact, already used in Section~\ref{sec.main} 
(see the proof of Lemma~\ref{2.2}):  
if a diagram $D$ 
on $m$ solid circles 
consists of  $m-1$ $\Y$'s and
has $W(D)\neq 0$, then 
at least one of the $m$ solid circles has exactly two univalent 
vertices on it. 

\begin{figure}[h]
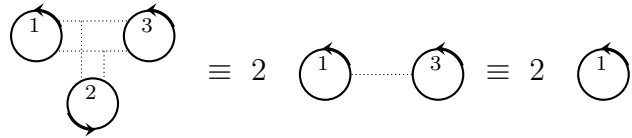
   
\begin{center}
\ExThree 
$\ \equiv  \ 2$ \ExThreei $\equiv  \ 2$ \ExThreeii

\caption{\label{FigExThree} Applying 
relation~(2c) of Figure~\ref{onevert}, 
we see that $W_4$ takes the value $2$ on this
diagram. Thus $\partial^2 F_3/\partial \mu_{123}^2=2$. By (ii), no
other diagram can contribute to $F_3$. This shows that
$F_3=\mu_{123}^2$.} 
\end{center}
\end{figure}

To check (i), see the computation in Figure~\ref{FigExThree}.

Finally,
formulas (iii) relate
the coefficients of monomials occurring in $F_m$
and the 
coefficients of monomials  
in $F_{m-2}$.  Note that these coefficients can be computed by
applying the Alexander-Conway weight system to the diagram
corresponding to a monomial. 
The identities in (iii) follow from
applying to 
this diagram 
relation~(2c) of Figure~\ref{onevert} (and smoothing the resulting
chords in the four remaining terms) 
in the case when the solid circle with the label $1$ has
exactly two univalent vertices on it.  
Let us elaborate. 
The operator 
$$
F_m \mapsto \mu_{1ij} \mu_{1kl} 
\left[\frac {\partial^2 \, F_m}{\partial \mu_{1ij}\, 
    \partial \mu_{1kl}}\right]_{v_1=0}
$$ 
projects $F_m$  onto 
the space of monomials of degree $m-1$ which contain
$\mu_{1ij}\mu_{1kl}$, but no other 
$\mu_{\alpha\beta\gamma}$ with $1\in\{\alpha,\beta,\gamma\}$. 
Such monomials correspond exactly to those diagrams where 
we {\em can apply} relation~(2c) 
with the distinguished solid circle in (2c) 
labelled by $1$ and the other four solid circles 
labelled by $i,j,k,l$.  
The operator 
$$
F_m \mapsto \left[\frac {\partial^2 \, F_m}{\partial \mu_{1ij}\,
    \partial \mu_{1kl}}\right]_{v_1=0}
$$
corresponds to 
{\em actually applying} (2c) to each of these diagrams,
and smoothing the chords in the resulting diagrams.
If in one of these 
new diagrams
the two endpoints of the chord lie on the same solid circle, then 
smoothing the chord increases the number of solid circles and 
the corresponding diagram
does not contribute (the number of remaining $\Y$'s is
too small).  If, however, the endpoints of the chord lie on the $i$th
and $j$th solid circle, with $i\neq j$, then the two circles are
replaced 
after the smoothing by their connected sum,
and such a term contributes a summand of the form
$$F_{m-2}(\ldots, v_i+v_j, \ldots)~.$$ 
There are three cases to consider 
corresponding
to whether 
the set $\{i,j,k,l\}$ has two, three, or four distinct elements. 
This leads to the three identities in~(iii). 

Finally, conditions 
 (i)---(iii) 
characterize $F_m$ uniquely, since
they allow to compute $F_m$ recursively. Indeed, (ii) shows that 
for every non-zero monomial in $F_m$ there is an index $p$ where one
can apply (iii)
and thus reduce the computation to $F_{m-2}$.

This completes the proof.\eproof

\begin{corollary}\label{cor7}
  The polynomial $F_m=F_m^{(2)}$ is equal to $\P_m^2$, the square of
  the Pfaffian-tree polynomial~(\ref{eq:pm}). 
\end{corollary}

\proof  This follows from Corollaries~6.5---6.7 of~\cite{MV},   
where it is shown that the polynomial $\P_m^2$ satisfies 
(and is determined by) the recursion relations in
Proposition~\ref{recursion}. \eproof

\begin{remark}
  {\em The proof that $\P_m^2$ satisfies the recursion relations is
  independent of the Pfaffian  Matrix-Tree Theorem for $3$-graphs of
  \cite{MV}. In fact, it follows from a three-term
  deletion-contraction relation for $\P_m$   which is proved using the
  interpretation of $\P_m$ as the spanning tree  generating
  function of  the complete $3$-graph with $m$ vertices.}   
\end{remark}

Thus we arrive at one of our main results (stated as 
Theorem~\ref{thm:alex-trees} in Section~\ref{sec:statmt}). 

\begin{corollary} Let $L$ be an $m$-component algebraically split
link, with triple Milnor linking numbers $\mu_{ijk}(L)$.  Let
$\nabla_L(z)=\sum_{i\geq 0} c_i(L) z^i$. 
Then $c_i(L)=0$ for $i\leq 2m-3$,  and 
$$c_{2m-2}(L)= \P_m(\mu_{ijk}(L))^2.$$ 
\end{corollary}

\proof Recall that the numbers $\mu_{ijk}(L)$ are the coefficients of
$\xi_2(L)$. Thus, the result follows from Proposition~\ref{5.3}, since
$F_m^{(2)}=\P_m^2$. \eproof

\subsection{Coefficients of $F_m^{(2)}$ and tree decompositions}
\

Let us fix an odd $m\ge 3$.
The  coefficients of the monomials 
in the polynomial $F_m^{(2)}=\P_m^2$ 
can be computed by counting decompositions of 
related 
$3$-graphs into pairs of spanning trees. 
Recall from Section~\ref{sec:statmt}  that   
edges   of a \emph{$3$-graph} have three (distinct) vertices and can
be visualized as   Y-shapes with the three vertices at their 
 endpoints (see~\cite{MV} for
details). A diagram $D$ whose components are $\Y$'s on $m$ ordered
solid circles defines a $3$-graph $G_D$  
with vertex set $\{1,2,\ldots ,m\}$ and edges 
given by these components.
Since the
edges of $3$-graphs
are not oriented, 
the assignment $D\mapsto G_D$ is not one-to-one. Note, however,
that the $3$-graph $G_D$ determines the diagram $D$ up to sign.  
Similarly, 
each  monomial $M$ in the variables $\mu_{ijk}$ gives
a $3$-graph $G_M$, and $G_M$ determines $M$ up to sign.

\begin{definition} \rm
Let $G$ be a $3$-graph. An {\em ordered  tree decomposition} of 
$G$  is a sub-$3$-graph  $T$ which is a tree and 
whose complement $T'$ is also a tree. 
\end{definition}

In the case when 
$G$ has $m$ vertices and $m-1$ edges, both trees in a tree
 decomposition must be spanning trees. 
In particular, for a diagram $D$ with  $m-1$\  $Y$-shaped components
or for a monomial $M$ of degree $m-1$ in the $\mu_{ijk}$'s each tree 
decomposition of the corresponding $3$-graph $G=G_D$ or $G=G_M$ 
gives two monomials $y_T$ and $y_{T'}$ of degree $(m-1)/2$
in the Pfaffian-tree polynomial $\P_m$.  
(Note that the monomials $y_T$ and $y_{T'}$ depend on $D$ or $M$ and
cannot be determined by the $3$-graph $G$ alone.)
Let $\varepsilon(T)$ and $\varepsilon(T')$ be the signs of these
monomials in $\P_m$. By counting each tree decomposition with the sign 
$$\varepsilon=\varepsilon(T)\varepsilon(T'),$$ 
we define the  \emph{algebraic number of ordered tree decompositions} 
of the diagram $D$ and of the monomial $M$.

With these definitions in hand,    the equality
   $F_m^{(2)}=\P_m^2$ (Corollary~\ref{cor7})
   and the definition of $\P_m$ as    the spanning tree generating function
   of the complete $3$-graph  $\Gamma_m$ (see Section~\ref{sec:statmt})  
immediately give the following  combinatorial rules for computing 
$W(D)$ and $F_m^{(2)}$. 
\begin{proposition}
\label{cor8}
\

(i) \ Let $D$ be a diagram consisting of $m-1$  $\Y$-shaped components
on $m$ solid circles. Then $W(D)$ is equal to
the   algebraic number of ordered tree decompositions of $D$.

(ii) \ Let $M$ be 
a monomial of degree $m-1$ in the  variables    $\mu_{ijk}$. 
Then the coefficient of $M$ in the polynomial $F_m^{(2)}$
is equal to the algebraic number of   ordered tree decompositions 
of $M$   divided by the symmetry factor $|Aut(G_M)|$,
where $Aut(G_M)$ is
the group of automorphisms of the $3$-graph $G_M$ 
inducing the identity map on the set of vertices of $G_M$. 
  \end{proposition}
\eproof

In our situation,
the cardinality  $|Aut(G_M)|$  is equal to $2^{d}$, 
where  $d$ is the number of $\mu_{ijk}$'s which occur twice in $M$ 
(up  to sign). In other words, $d$ is the number of (unordered) triples of 
vertices in $G_M$ with $2$ edges  attached to them.
(Note that if $M$ has non-zero coefficient in $F_m^{(2)}=\P_m^2$,
they $G_M$ can have at most two edges with the same vertex set.)

To see 
the need for
the symmetry factor  $|Aut(G_M)|$  
in (ii), consider for example the monomial $y_{123}^2$.  
It occurs
in $\P_3^2$ with coefficient $1$, 
but the weight system $W$ takes the value
$2$ on the
diagram in Figure~\ref{FigExThree}
 which is dual to $y_{123}^2$. This corresponds to the fact that the 
associated $3$-graph $G_M$ has two ordered tree decompositions and  
$|Aut(G_M)|=2$.

Here are two more examples to illustrate this. The
monomial $$M=y_{123}^2y_{245}y_{345}$$ has four ordered tree
decompositions, each 
contributing $+1$, 
and $|Aut(G_M)|=2$. 
Therefore, $M$ appears  
in $\P_3^2$ with  coefficient $2$. 
Finally, the monomial 
$$ M=y_{145}\, y_{146}\, y_{256}\, y_{257}\, y_{347}\, y_{367}$$ 
has six  ordered tree decompositions (again each contributing $+1$) and
$|Aut(G_M)|=1$; therefore it has coefficient  $6$ in $\P_7^2$.
\medskip

Let us finish this subsection by illustrating how to compute the
coefficients of monomials in $\P_m^2$ using the recursion relations of
Proposition~\ref{recursion}. 
 For example, let us compute the monomials  in $F_5$ which contain
    $\mu_{123}\mu_{145}$, but no other $\mu_{ijk}$  
with $1\in\{i,j,k\}$. The coefficient of such a monomial in $F_5$ is
    equal to its coefficient in 
 \begin{equation}
      \label{eq:example}
\mu_{123}\,\mu_{145}\, \left[\frac {\partial^2 \, F_5}{\partial
    \mu_{123}\, \partial \mu_{145}}\right]_{v_1=0}~.
 \end{equation}
 Applying~\ref{recursion}(iii), 
we see that~(\ref{eq:example}) is equal to
\begin{align*} &\mu_{123}\,\mu_{145}\, \bigl( F_{3}( v_3+v_4,v_2, v_5) +
          F_{3}( v_2 +v_5, v_3, v_4) \\
&\hskip 4.5cm -F_{3}( v_2+v_4,v_3, v_5) - 
F_{3}( v_3 +v_5, v_2, v_4)\big)\\ 
&=\mu_{123}\,\mu_{145}\,\big((\mu_{325}+\mu_{425})^2
+(\mu_{234}+\mu_{534})^2 \\
&\hskip 6cm - (\mu_{235}+\mu_{435})^2 -(\mu_{324}+\mu_{524})^2\big)\\
&= 2\, \mu_{123}\,\mu_{145}\,(\mu_{234} +\mu_{235})(\mu_{245}+\mu_{345})~.
\end{align*}

\begin{figure}[h]
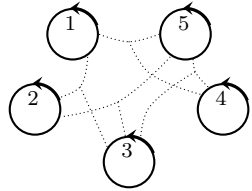
   
\begin{center}
\ExFive 

\caption{\label{FigExFive} 
This diagram contributes the term 
$2\, \mu_{123}\,\mu_{145} \,\mu_{235}\,\mu_{345}$ to $F_5=\P_5^2$. 
The associated $3$-graph can be decomposed into spanning trees corresponding to
$\mu_{123}\,\mu_{345}$ and $\mu_{145}\,\mu_{235}$.} 
\end{center}
\end{figure}

For example, we see that the coefficient of 
$\mu_{123}\,\mu_{145} \,\mu_{235}\,\mu_{345}$
     in $F_5$ is equal to $2$ (see Figure \ref{FigExFive}).

\subsection{The case $n\geq 3$.}

\

Denote
by $\overline{C}^t_2(m)$ the quotient ${C}^t_2(m)/W_0$ where
$W_0$ is the subspace of ${C}^t_2(m)=\Lambda^3V$ defined as
follows. Let $\Y_{ijk}$  be the basis of $\Lambda^3V$ dual to
$\mu_{ijk}\in \Lambda^3V^*$. Then $W_0$ is the subspace  generated by 
elements  of the form
$$
(\Y_{ijk} -\Y_{ijl}) - (\Y_{jkl}-\Y_{ikl})$$ 
for every set $\{i,j,k,l\}$
of four distinct vertices (see Figure \ref{figmu}).

\begin{figure}[h]   
\begin{center}
\MUa\  $-$ \ \MUb \ = \ \MUc \ $-$ \ \MUd

\caption{\label{figmu} The relations in $\overline{C}^t_2(m)$.}
\end{center}
\end{figure}
The Pfaffian-tree polynomial $\P_m$ can be viewed as a polynomial function 
on $C_2^t(m)$.
It is shown in \cite[Proposition 6.10]{MV} that 
$\P_m$ descends to a well-defined polynomial on 
$\overline{C}^t_2(m)$. Hence, the same is true for 
$F_m^{(2)}=\P_m^2$. Alternatively, the fact that $F_m^{(2)}$ is well-defined
on $\overline{C}^t_2(m)$ 
can be proved
from the definition of $F_m^{(2)}$
via the Alexander-Conway weight system and the relations in
Figure~\ref{Tree4}.

\begin{proposition} \label{ACred} The Alexander-Conway relations 
define a linear map   
\begin{equation}
  \label{eq:phi}
\phi_n : C^t_n(m) \rightarrow 
\begin{cases} C_1(m)\ \ \ \text{if $n$ is odd}\\
\overline{C}^t_2(m)\ \ \ \text{if $n$ is even}
\end{cases}
\end{equation}
 such that 
\begin{equation}\label{Fm}
F_m^{(n)}=\begin{cases} F_m^{(1)} \circ \phi_n =\ \D_m \circ \phi_n~,
 \ \ \   \text{if $n$ is odd}\\ 
F_m^{(2)} \circ \phi_n= \ \P_m^2 \circ \phi_n~,
\ \ \ \text{if $n$ is even}
\end{cases}
\end{equation}
\end{proposition}
\proof Let $D$ be a labelled tree diagram of degree $n$. We define
  $\phi_n(D)$   to be the unique element of  $C_1(m)$ or 
$\overline{C}^t_2(m)$, depending on whether $n$ is odd or even, 
such that $D\equiv \phi_n(D)$. 
The existence of $\phi_n(D)$ was proved in 
the tree reduction Lemma~\ref{lemred}.  

To prove uniqueness of $\phi_n(D)$, we need to show that applying the 
tree reduction relations (see Equation~(\ref{Relconseq}) and 
Figures~\ref{Tree4}   and~\ref{reducetrees}) in any order always gives
the same result. This verification is an easy application of the
diamond lemma. Note that the apparent indeterminacy coming from the
relation in Figure~\ref{Tree4} is compensated by taking the 
quotient $\overline{C}^t_2(m)={C}^t_2(m)/W_0$ in  
the case when $n$ is even.
Finally, Equation~(\ref{Fm}) follows immediately from the definition
of $\phi_n$.   This completes the proof.\eproof

\medskip

Thus, we arrive at the general formula of 
Theorem~\ref{thm:general}.

\begin{corollary} \label{r7.7} Let $L$ be an $m$-component link 
with vanishing Milnor invariants of degree $\le n-1$ 
(\emph{i.e.}\ $L\in \L_n(m)$). 
Then $c_i(L)=0$ for $i< n(m-1)$, and 
$$c_{n(m-1)}(L)=
\begin{cases} 
(\D_m   \circ \phi_n)(\xi_n(L)) &=
\D_m(\ell^{(n)}_{ij}(L))~,  \qquad \text{if $n$ is odd}\\ 
(\P_m \circ \phi_n)(\xi_n(L))^2 &=
(\P_m(\mu^{(n)}_{ijk}(L)))^2~,  \ \ \text{if $n$ is even}
\end{cases}$$
where $\phi_n$ is the reduction map~(\ref{eq:phi})
and $\xi_n(L)$ is the universal degree-$n$ Milnor invariant~(\ref{eq:xi}) 
of $L$. 
\end{corollary}
Here we used the linear forms on $C_n^t(m)$ given by
$$
    \ell^{(n)}_{ij}=\phi_n\circ \ell_{ij} 
$$ 
(when $n$ is odd)
and 
    $$\mu^{(n)}_{ijk}=\mu_{ijk}\circ\tilde{\phi}_n,
    $$
where 
$$
   \tilde{\phi}_n: C_n^t(m)\to C_2^t(m)
$$ 
is an arbitrary lift of the map $\phi_n$ to $C_2^t(m)=\Lambda^3(V)$
(when $n$ is even).\footnote{The linear form $\mu^{(n)}_{ijk}$ depends 
on the lift $\tilde{\phi}_n$, but $\P_m(\mu^{(n)}_{ijk})$ is independent
of this choice.}

\proof This follows from Propositions~\ref{5.3} and~\ref{ACred}. \eproof

\section{The case when both $n$ and $m$ are even} 
\label{sec:even}

If both $n$ and $m$ are even, then by~(\ref{eq:parvanish})
the coefficient $c_{n(m-1)}(L)$ of 
the Alexander-Conway polynomial 
$\nabla_L$ is zero for every $m$-component link $L$. 
Therefore, by Corollary~\ref{r7.7},
if $L\in \L_n(m)$, then
the first non-vanishing coefficient of $\nabla_L$ can 
occur only in degree $n(m-1)+1$.
The Vanishing Lemma~\ref{2.1} and the methods of
Section~\ref{sec:coeff} give an expression of this coefficients as a
polynomial in Milnor numbers,
analogous to the formula 
 $c_{n(m-1)}(L)=F_m^{(n)}(\xi_n(L))$ of Proposition~\ref{5.3}.   
\begin{proposition}
There exists
a polynomial on $C_n^t(m) \oplus C_{n+1}^t(m)$
\begin{equation}\label{GMN2}
        G_m^{(n)}\in  (S^{m-2}C_n^t(m)^*) \otimes  C_{n+1}^t(m)^*
\subset S^{m-1}(C_n^t(m)^* \oplus C_{n+1}^t(m)^* ),
\end{equation} 
such that
\begin{equation}\label{Gmn}
c_{n(m-1)+1}(L)=G_m^{(n)}(\xi_n(L),\xi_{n+1}(\s)).
\end{equation}
       Here
$\xi_n(L)$ is the universal 
degree-$n$ Milnor invariant of $L$ 
and $\xi_{n+1}(\s)\in C_{n+1}^t(m)$ corresponds to
$Z_{n+1}^t(\s)\in \A^t_{n+1}(m)$, where $\sigma$ is any string link
representative of $L$.
\end{proposition}

\proof
The homogeneous polynomial  $G_m^{(n)}$ is defined
similarly to the polynomial $F_m^{(n)}$
by lifting a disjoint union of $m-2$ labelled diagrams from $C_n^t(m)$ 
and one diagram from $C_{n+1}^t(m)$ to  $\A_{n(m-1)+1}(\amalg_m S^1)$
and applying $W$. 
Independence of the result on the choice of
lift follows from the STU 
relations and the Vanishing Lemma similarly to the proof of
Lemma~\ref{Fdef}. 
\eproof

Although $\xi_{n+1}(\s)\in C_{n+1}^t(m)$ may depend on the string link
representative $\s$ and not just on the link $L$, formula (\ref{Gmn})
shows that $G_m^{(n)}(\xi_n(L),\xi_{n+1}(\s))$  depends only on $L$, 
since it is equal to the coefficient $c_{n(m-1)+1}(L)$ of $\nabla_L$. 

For example, if $L$ is an algebraically split $m$-component link 
where $m$ is even, then the first non-vanishing coefficient
$c_{2m-1}(L)$ is a polynomial of degree $m-2$ in the triple linking
numbers $\mu_{ijk}(L)$ (which are the coefficients of $\xi_2(L)$) and
of degree one in the quadruple linking numbers $\mu_{ijkl}(\s)$ 
(which are encoded by the coefficients of $\xi_3(\s)$). 

\begin{figure}[h]
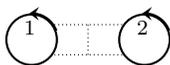
   
\begin{center}
\ExFour

\caption{\label{figmu1122} The only diagram contributing to
$G_2^{(2)}$.} 
\end{center}
\end{figure}

If $m=2$, there is only one possible diagram contributing to
$G_2^{(2)}$ (see Figure~\ref{figmu1122}). The coefficient of this
diagram in $Z^t_3(L)$ is $\pm \frac{1}{2}\mu_{1122}(L)$ (see
\cite[section 8]{HM}). Since $W$ takes the value $-2$ on this diagram,
we recover the well-known fact \cite{Co} that for $2$-component
algebraically split links $$c_3(L)=\pm \mu_{1122}(L)$$ (the sign 
depends on conventions which we don't want to specify here).

\begin{figure}[h]
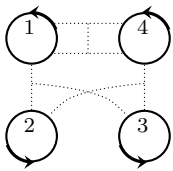
   
\begin{center}
\ExFouri

\caption{\label{figFour} A diagram contributing to $G_4^{(2)}$.}
\end{center}
\end{figure}

The existence of some expression for $c_{n(m-1)+1}(L)$ as a polynomial
in Milnor numbers also  follows from Levine's general determinantal formula
in \cite{Le2} (see also Traldi~\cite{Tr2}). In fact, it is even possible to
infer from Levine's  formula that the polynomial $G_m^{(n)}$ must be of 
the form given in~(\ref{GMN2}), \emph{i.e.}\ it is homogeneous of degree 
$m-2$ in
Milnor numbers of length $n+1$ 
and of degree one in 
Milnor numbers
of length $n+2$.
But while Levine's formula  describes 
the polynomial  $G_m^{(n)}$ as a determinant,  our methods again lead 
to an
expression for it  in terms of tree
decompositions. 

Indeed, the computation of the coefficients of
$G_m^{(n)}$ can be reduced to those of $F_m^{(n)}$.
Consider, for example,
the case $n=2$. 
Then we only need to look at diagrams like the one 
in Figure~\ref{figFour}.
Applying relation~\eqref{Relb} to the
 vertical  dashed edge of the $H$-shaped component of this diagram, 
we may insert an additional  solid circle into that edge 
without changing the value of the Alexander-Conway weight system. Thus, 
we see that the diagram in
Figure~\ref{figFour} is equivalent to a diagram
consisting of four $\Y$'s on five solid circles. 
Similarly, the diagram in
Figure~\ref{figmu1122} is equivalent to a diagram
consisting of two $\Y$'s on three solid circles  
(which is equal to minus the diagram in Figure~\ref{FigExThree}). 

The same argument shows that in general, a diagram consisting of $m-2$ $\Y$'s
 and one $H$-shaped  component  on $m$ 
solid circles is equivalent to a diagram consisting of $m$ $\Y$'s on
$m+1$ solid circles,
and thus, we arrive at the following result.

\begin{proposition}
\label{cor9} Let $D$ be a labelled diagram (with labels
   from the set $1,2,\ldots, m$)  which has $m-2$ \ $\Y$-shaped  components
and one $H$-shaped   component ${}^\ell_iH^k_j$ (with obvious notation).  
Then $$W(D)=W(D'),$$
 where $D'$ is the diagram (with an additional label 
$0$) obtained from $D$
 by replacing the component ${}^l_iH^k_j$ with two $Y$-shaped
components  $Y_{0jk}$ and  $Y_{0\ell i}$.
\end{proposition}
\eproof

It follows that $W(D)$ is 
equal to
the algebraic number of tree decompositions of the  
diagram
$D'$ (see 
Proposition~\ref{cor8}(i)). 
We leave it to the reader to translate this fact
into a formula for the polynomial $G_m^{(2)}$ in terms of monomials involving
$m-2$ $\mu_{\alpha\beta\gamma}$'s and one ${}^\ell_i\eta^k_j$, where 
${}^\ell_i\eta^k_j$ is the dual  of ${}^\ell_iH^k_j$.
For example, one has 
$$G_2^{(2)}= - 2 \ {}^1_2\eta^1_2~.$$
To do this for $m\geq 4$,
one has to take into account the symmetry factors as in 
Proposition~\ref{cor8}(ii).
A certain technical point arises here because of the fact
that the dual of a diagram does  
not make sense in $C_n^t(m)^*$ for $n\ge 3$. 
(This point did not arise in
the description of the
polynomial $F_m^{(2)}=\P_m^2$, since the IHX relation is not needed to 
describe  $\Y$-shaped labelled diagrams.) The 
standard solution is to consider the space  $\tilde{C}_n^t(m)$
  of labelled tree diagrams modulo AS relations (but not
modulo IHX relations).  
The dual ${}^\ell_i\eta^k_j$ of ${}^\ell_iH^k_j$ 
makes sense in      $\tilde{C}_3^t(m)^*$,
 and one obtains from Proposition~\ref{cor9} 
a tree decomposition-type formula
    for the coefficients of monomials of the polynomial $G_m^{(2)}$.
Each of these monomials 
is a function on
        $\tilde{C}_2^t(m) \oplus  \tilde{C}_{3}^t(m)$
and in general does  
not descend to   $C_2^t(m) \oplus  C_{3}^t(m)$.
       However, 
their sum is equal to the polynomial $G_m^{(2)}$   which satisfies 
the IHX relations and  
therefore is a well-defined function
on the quotient.

\end{document}